\documentclass{amsart}
\usepackage{amsfonts}
\usepackage{pstricks,pstcol,pst-plot,pst-node,pst-tree,amssymb}
\usepackage{pstricks-add}
\usepackage{graphicx}
\usepackage{caption}
\usepackage{mathrsfs}

\theoremstyle{definition}

\theoremstyle{definition}

\theoremstyle{definition}

\theoremstyle{definition}

\theoremstyle{definition}

\theoremstyle{definition}

\numberwithin{equation}{section}

\input{tcilatex}

\begin{document}
\title[Pansu spheres]{UMBILICITY AND CHARACTERIZATION OF PANSU SPHERES IN
THE HEISENBERG GROUP}
\author{Jih-Hsin Cheng}
\address{Institute of Mathematics, Academia Sinica, Taipei and National
Center for Theoretical Sciences, Taipei Office, Taiwan, R.O.C.}
\email{cheng@math.sinica.edu.tw}
\thanks{}
\author{Hung-Lin Chiu}
\address{Department of Mathematics, National Central University, Chung Li,
32054, Taiwan, R.O.C.}
\email{hlchiu@math.ncu.edu.tw}
\urladdr{}
\author{Jenn-Fang Hwang}
\address{Institute of Mathematics, Academia Sinica, Taipei, Taiwan, R.O.C.}
\email{majfh@math.sinica.edu.tw}
\thanks{}
\author{Paul Yang}
\address{Department of Mathematics, Princeton University, Princeton, NJ
08544, U.S.A.}
\email{yang@Math.Princeton.EDU}

\begin{abstract}
For $n\geq 2$ we define a notion of umbilicity for hypersurfaces in the
Heisenberg group $H_{n}$. We classify umbilic hypersurfaces in some cases,
and prove that Pansu spheres are the only umbilic spheres with positive
constant $p($or horizontal)-mean curvature in $H_{n}$ up to Heisenberg
translations.
\end{abstract}

\maketitle



\section{Introduction and statement of the results}

In classical differential geometry, we have the notion of umbilicity for a
point in a hypersurface of the Euclidean space $R^{n}$. A connected, closed
umbilic hypersurface of $R^{n}$ (i.e., all the points are umbilic) is shown
to be a sphere. On the other hand, we have the Alexandrov theorem which says
that a closed (compact with no boundary) hypersurface of positive constant
mean curvature in $R^{n}$ must be a sphere. The original proof of
Alexandrov's theorem (\cite{Al}) is based on a reflection principle. Reflect
the hypersurface $S$ across a hyperplane $P.$ Move $P$ until the reflected
hypersurface touches the original hypersurface $S.$ The reflected
hypersurface must coincide with $S$ by the strong maximum principle.
Analytic proofs of Alexandrov's theorem were given much later. In 1991
Montiel and Ros (\cite{MR}) gave a relatively elementary proof through the
characterization of spheres by the umbilicity.

For a hypersurface in the Heisenberg group $H_{n}$ (see Section 2 for some
basic material)$,$ we can still talk about mean curvature, called $p($or
horizontal)-mean curvature $H$ (see Section 2 for the definition). A
hypersurface defined by such $H$ $=$ $0$ is called $p$(horizontal)-minimal.
Such $p$-minimal hypersurfaces or hypersurfaces with prescribed $p$-mean
curvature have been extensively studied in the last ten years (see, for
instance, \cite{Pau04}, \cite{CH04}, \cite{CHMY04}, \cite{Pauls}, \cite{CHY1}%
, \cite{Ri}, \cite{RR}, \cite{CHY2}, \cite{CH10}, \cite{CHMY12}, \cite{Rit}, 
\cite{CH14}, and references therein).

By analogy with the Euclidean situation, we can ask if an Alexandrov-type
theorem holds for the Heisenberg situation. The reflection principle doesn't
seem to work generally in this situation. In the case $n$ $=$ $1,$ Ritore
and Rosales (\cite{RR}) showed that an Alexandrov-type theorem still holds.
Their proof relies on the analysis of characteristic curves and singular set
developed in \cite{CHMY04}. For $n$ $\geq $ $2,$ on the other hand, we may
invoke the method of Montiel and Ros to study the Alexandrov-type problem.
So the first thing is to characterize, in this case, Pansu spheres (having
positive constant $p$-mean curvature; see (\ref{1.6})) in terms of some
notion of umbilicity. In this paper, we give a definition of umbilicity. We
classify umbilic hypersurfaces in some cases, and carry out a
characterization of Pansu spheres in $H_{n}$.

Let $\Sigma $ be a $C^{2}$ smooth (further assume the regular part is $%
C^{\infty }$ smooth; see below) hypersurface of the Heisenberg group $H_{n}$%
. Throughout this paper, we always assume $\Sigma $ is immersed and $n$ $%
\geq $ $2$. Let $\xi $ ($J$, resp.) denote the standard contact ($CR$,
resp.) structure on $H_{n},$ defined by the kernel of the contact form 
\begin{equation*}
\Theta =dt+\sum_{j=1}^{n}(x_{j}dy_{j}-y_{j}dx_{j})
\end{equation*}%
\noindent (see \cite{CL}, \cite{Lee}, or Section 2). A point $p$ $\in $ $%
\Sigma $ is called singular if $\xi $ $=$ $T\Sigma $ at $p$. Otherwise $p$
is called regular or nonsingular (i.e., $\xi $ is transversal to $T\Sigma )$%
. Let $S_{\Sigma }$ denote the set of singular points, which is a closed
subset of $\Sigma .$ We will further assume the regular part $\Sigma
\backslash S_{\Sigma }$ is $C^{\infty }$ smooth. For a regular point, we
define $\xi ^{\prime }$ $\subset $ $\xi $ $\cap $ $T\Sigma $ by%
\begin{equation}
\xi ^{\prime }=(\xi \cap T\Sigma )\cap J(\xi \cap T\Sigma ).  \label{1.1}
\end{equation}

\noindent Let ($\xi ^{\prime }$)$^{\bot }$ denote the space of vectors in $%
\xi $, perpendicular to $\xi ^{\prime }$ with respect to the Levi metric $G$
:$=$ $\frac{1}{2}d\Theta (\cdot ,J\cdot )$ $=$ $%
\sum_{j=1}^{n}[(dx_{j})^{2}+(dy_{j})^{2}].$ It is not hard to see $\dim (\xi 
$ $\cap $ $T\Sigma )$ $\cap $ ($\xi ^{\prime }$)$^{\bot }$ $=$ $1.$ Take $%
e_{n}$ $\in $ ($\xi $ $\cap $ $T\Sigma )$ $\cap $ ($\xi ^{\prime }$)$^{\bot
} $ of unit length. Define the horizontal normal $e_{2n}$ :$=$ $Je_{n}$. Let 
$\nabla $ denote the pseudohermitian connection associated to $(J,\Theta )$
(see Section 2 for an explanation). Observe that $\nabla _{e_{n}}e_{2n}$ $%
\in $ $\xi $ is perpendicular to $e_{2n}.$ So we can write $-\nabla
_{e_{n}}e_{2n}$ $=$ $le_{n}$ modulo $\xi ^{\prime }$ for some function $l.$
Now define the vector field $X_{n}$ $\in $ $\xi ^{\prime }$ by%
\begin{equation}
X_{n}:=\nabla _{e_{n}}e_{2n}+le_{n}.  \label{1.2}
\end{equation}%
\noindent This vector field is uniquely defined on the regular part of $%
\Sigma $. Note that if $p\in \Sigma $ is a regular point such that $%
X_{n}(p)=0$, then we have 
\begin{equation}
\left( -\nabla e_{2n}+\alpha J^{\prime }\right) (\xi ^{\prime })\subset \xi
^{\prime },
\end{equation}%
\noindent (see Proposition 2.3) where 
\begin{equation}
J^{\prime }:=J\text{ on }\xi ^{\prime }\text{ and }J^{\prime }e_{n}:=0
\end{equation}%
\noindent (cf. (\ref{A5.1})). Hence we can regard this operator $-\nabla
e_{2n}+\alpha J^{\prime }$ originally defined on $\xi \cap T\Sigma $ (see (%
\ref{A5.2})) as an endomorphism on $\xi ^{\prime }$. This symmetric second
fundamental form or shape operator first appeared in Ritor\'{e}'s paper (see
page 52 in \cite{Rit}). Conversely, if $\xi ^{\prime }$ is invariant under
the operator $-\nabla e_{2n}+\alpha J^{\prime }$, then $X_{n}=0$ (see also
Proposition 2.3). In addition, it is self-adjoint (see Proposition 2.2). So
we immediately have the following result.

\bigskip

\textbf{Proposition 1.1.} Let $p$ be a regular point of $\Sigma $ such that $%
X_{n}(p)=0$. There are scalars 
\begin{equation*}
\lambda _{\beta },\ \lambda _{n+\beta },\ \ \ 1\leq \beta \leq n-1
\end{equation*}%
\noindent and an orthonormal basis 
\begin{equation*}
e_{\beta },\ e_{n+\beta },\ \ \ 1\leq \beta \leq n-1
\end{equation*}%
\noindent of $\xi ^{\prime }(p)$ such that 
\begin{equation}
\left( -\nabla e_{2n}+\alpha J^{\prime }\right) (e_{j})=\lambda _{j}e_{j},\
\ \text{for}\ \ 1\leq j\leq 2n-1,\ j\neq n.
\end{equation}

\bigskip

\textbf{Definition 1.2}. A regular point $p\in \Sigma $ is called an \textbf{%
umbilic} point if

(1) $\left( -\nabla e_{2n}+\alpha J^{\prime }\right) (\xi ^{\prime })\subset
\xi ^{\prime }$, and

(2) $\lambda _{1}=\cdots =\lambda _{n-1}=\lambda _{n+1}=\cdots =\lambda
_{2n-1}$.

\bigskip

If all regular points of $\Sigma $ are umbilic, we call $\Sigma $ an \textbf{%
umbilic hypersurface} of the Heisenberg group $H_{n}$. We often use $\lambda 
$ (or $k$) to denote the common eigenvalue in (2) of Definition 1.2.

For any $\lambda $\TEXTsymbol{>}$0$, the Pansu sphere $S_{\lambda }$ is the
union of the graphs of the functions $f$ and $-f$, where%
\begin{equation}
f(z)=\frac{1}{2\lambda ^{2}}\left( \lambda |z|\sqrt{1-\lambda ^{2}|z|^{2}}%
+\cos ^{-1}{\lambda |z|}\right) ,\ \ \ |z|\leq \frac{1}{\lambda }.
\label{1.6}
\end{equation}

\noindent It~is~known~that $S_{\lambda }$ has $p$-(or horizontal) mean
curvature $H$ $=$ $2n\lambda $ (see Section 2 for basic definitions and
Example 3.2 for more discussion; also see, for instance, \cite{Rit}). We say
that $\Sigma $ is congruent with a Pansu sphere if after a Heisenberg
translation, $\Sigma $ coincides with $S_{\lambda }$ for some $\lambda $ $>$ 
$0.$

\bigskip

\textbf{Theorem A}. \textit{Suppose }$\Sigma $\textit{\ is a closed,
connected umbilic hypersurface of }$H_{n}$\textit{\ (}$\mathit{n\geq 2)}$ 
\textit{with positive constant }$p$\textit{-mean curvature and nonvanishing
Euler number. Then }$\Sigma $\textit{\ is congruent with a Pansu sphere.}

\bigskip

\textbf{Corollary} $A^{\prime }.$ \textit{Suppose }$\Sigma $\textit{\ is
homeomorphic to the sphere }$S^{2n}.$ \textit{Suppose }$\Sigma $\textit{\ is}
\textit{an umbilic hypersurface of }$H_{n}$\textit{\ with positive constant }%
$p$\textit{-mean curvature. Then }$\Sigma $\textit{\ is congruent with a
Pansu sphere.}

\bigskip

Note that $S^{2n}$ is closed, connected, and having nonzero Euler number. So
Corollary $A^{\prime }$ follows from Theorem A immediately.

\bigskip

\textbf{Theorem 1.3.} \textit{Suppose }$\Sigma $\textit{\ is a closed,
connected umbilic hypersurface with }$l=2k$\textit{. Then }$\Sigma $\textit{%
\ is congruent with a Pansu sphere }$S_{\lambda }$\textit{\ with }$\lambda
=k $\textit{.}

\bigskip

\textbf{Lemma B}. \textit{Suppose }$\Sigma $\textit{\ is a connected umbilic
hypersurface of }$H_{n}$ \textit{with positive constant }$p$\textit{-mean
curvature, containing a singular point. Then }$l$\textit{\ }$=$\textit{\ }$%
2k $\textit{\ on }$\Sigma \backslash S_{\Sigma }.$

\bigskip

\textbf{Theorem 1.4.} \textit{Suppose }$\Sigma $\textit{\ is an umbilic
hypersurface with }$l=2k$\textit{. Then }$k$\textit{, and hence }$l$\textit{%
, are constants on the whole regular part of }$\Sigma $\textit{. Moreover,
if }$\Sigma $\textit{\ is connected and there exists a singular point }$p\in
\Sigma $\textit{, then }$\Sigma $\textit{\ is either congruent with part of
a Pansu sphere or congruent with part of a hyperplane orthogonal to the }$t$%
\textit{-axis.}

\bigskip

Theorem 1.3 is an immediate consequence of Theorem 1.4. This is because that
If $\Sigma $ is closed, then it must contain a singular point. Otherwise
Proposition 4.5 would imply that $\Sigma $ is foliated by geodesics, a
contradiction to compactness of $\Sigma $. Also the constant $l$ must be
positive. On the other hand, Proposition 4.1 shows that this singular point
is isolated, hence $\Sigma $ is congruent with a Pansu sphere $S_{\lambda }$
with $\lambda =k$. It was shown in \cite{LM} that for a rotationally
invariant hypersurface in $H_{n}$ with $l$ $=$ $2k$ we have the same
conclusion as in Theorem 1.4. Note that rotationally invariance implies
umbilicity by Proposition 3.1.

In Example 3.4, we introduce two kind of umbilic hypersurfaces with $\alpha
=0$. The hypersurface $\Sigma _{S^{2n-1}(c)}$ satisfies $l=k=\frac{1}{c}$.
The other one $\Sigma _{E}$ satisfies $k=l=0$. Conversely, we have the
following result.

\bigskip

\textbf{Theorem 1.5}. \textit{Suppose }$\Sigma $\textit{\ is an umbilic
hypersurface with }$\alpha =0$\textit{. Then }$k$\textit{\ is a constant on }%
$\Sigma $\textit{. Moreover, if }$\Sigma $\textit{\ is connected and }$k>0$%
\textit{, then }$l=k$\textit{, and hence }$\Sigma $\textit{\ is congruent
with part of the hypersurface }$\Sigma _{S^{2n-1}(c)}$\textit{\ with }$c=%
\frac{1}{k}$\textit{. If }$\Sigma $\textit{\ is connected and }$k=l=0$%
\textit{, then }$\Sigma $\textit{\ is congruent with part of the
hypersurface }$\Sigma _{E}$\textit{\ for some hyperplane }$E$\textit{.}

\bigskip

In Section 2 we give a sketch of the basic theory of hypersurfaces in $%
H_{n}. $ In particular, we discuss the symmetry property of the second
fundamental form. We end up defining a symmetric second fundamental form or
shape operator. In Section 3 we show that rotationally invariance implies
umbilicity and give examples including Pansu spheres, Heisenberg spheres,
and umbilic hypersurfaces with $\alpha $ $=$ $0.$

In Section 4 we study important properties of umbilic hypersurfaces and
prove Theorem 1.4 and Theorem 1.5. We postpone the proof of Proposition 4.2
to Section 5. Included in Proposition 4.2 are many useful formulas for
umbilic hypersurfaces. In Section 6 we study an ODE system associated to an
umbilic hypersurface. A complete understanding of this ODE system (Lemma
6.1) helps us to give a proof of Lemma B. We can finally prove Theorem A in
Section 7. Besides, we observe examples of Sobolev extremals whose level
sets are umbilic hypersurfaces and pose a question whether each level set of
a Sobolev extremal is umbilic.

\bigskip

\textbf{Acknowledgments}. J.-H. C. (P. Y., resp.) is grateful to Princeton
University (Academia Sinica in Taiwan, resp.) for the kind hospitality.
J.-H. C., H.-L. C., and J.-F. H. would like to thank the Ministry of Science
and Technology of Taiwan for the support of the following research projects:
NSC 101-2115-M-001-015-MY3, NSC 100-2628-M-008-001-MY4, and NSC
102-2115-M-001-003-MY2, resp.. P. Y. would like to thank the NSF of the
United States for the grant DMS-1104536. We thank the referee for careful
reading of the argument and pointing out a small gap in the previous
version.We would also like to thank Ms.Yu-Tuan Lin for the computer
assistance to draw Figure 6.1 and Figure 6.2.

\bigskip

\section{Basic theory of hypersurfaces in $H_{n}$}

The Heisenberg group $H_{n}$ is $R^{2n+1},$ as a set, together with the
group multiplication%
\begin{eqnarray*}
&&(x_{1},..,x_{n},y_{1},..,y_{n},t)\circ (\tilde{x}_{1},..,\tilde{x}_{n},%
\tilde{y}_{1},..,\tilde{y}_{n},\tilde{t}) \\
&=&(x_{1}+\tilde{x}_{1},..,x_{n}+\tilde{x}_{n},y_{1}+\tilde{y}_{1},..,y_{n}+%
\tilde{y}_{n},t+\tilde{t}+\sum_{j=1}^{n}(y_{j}\tilde{x}_{j}-x_{j}\tilde{y}%
_{j})).
\end{eqnarray*}

\noindent $H_{n}$ is a $(2n+1)$-dimensional Lie group. Any left invariant
vector field is a linear combination of the following basic vector fields:%
\begin{eqnarray*}
\mathring{e}_{j} &=&\frac{\partial }{\partial x_{j}}+y_{j}\frac{\partial }{%
\partial t},\mathring{e}_{n+j}=\frac{\partial }{\partial y_{j}}-x_{j}\frac{%
\partial }{\partial t},1\leq j\leq n \\
\text{and }T &=&\frac{\partial }{\partial t}.
\end{eqnarray*}

\noindent The standard contact structure $\xi $ on $H_{n}$ is the subbundle
of $TH_{n},$ spanned by $\mathring{e}_{j}$ and $\mathring{e}_{n+j},$ $1\leq
j\leq n.$ Or equivalently we can define $\xi $ to be the kernel of the
standard contact form%
\begin{equation*}
\Theta =dt+\sum_{j=1}^{n}(x_{j}dy_{j}-y_{j}dx_{j}).
\end{equation*}

\noindent The standard $CR$ structure on $H_{n}$ is the almost complex
structure $J$ defined on $\xi $ by%
\begin{equation*}
J(\mathring{e}_{j})=\mathring{e}_{n+j}\text{ and }J(\mathring{e}_{n+j})=-%
\mathring{e}_{j}.
\end{equation*}

Recall the pseudohermitian structure $(J,\Theta )$ on $H_{n}$ (\cite{Web}, 
\cite{Lee}) as follows$.$ Let $\nabla $ denote the pseudohermitian
connection. It has the following good property:%
\begin{equation*}
\nabla \mathring{e}_{j}=\nabla \mathring{e}_{n+j}=\nabla T=0
\end{equation*}%
\noindent for $1$ $\leq $ $j$ $\leq $ $n.$ Write $\xi \otimes C$ $=$ $%
T_{1,0}\oplus T_{0,1}$ where $T_{1,0}$ ($T_{0,1},$ resp.) is the eigenspace
of $J$ with eigenvalue $i$ ($-i,$ resp.) (at each point). Then there exist
complex-valued 1-forms (called unitary coframe) $\theta ^{\beta },$ $1\leq
\beta \leq n,$ which annihilate $T_{0,1}$ and $T,$ such that%
\begin{equation}
d\Theta =i\sum_{\beta =1}^{n}\theta ^{\beta }\wedge \theta ^{\bar{\beta}}
\label{A1}
\end{equation}

\noindent ($\theta ^{\bar{\beta}}$ means the complex conjugate of $\theta
^{\beta }).$ Let $\theta _{\beta }$ $^{\gamma }$ denote the pseudohermitian
connection forms such that 
\begin{eqnarray}
d\theta ^{\beta } &=&\theta ^{\gamma }\wedge \theta _{\gamma }\text{ }%
^{\beta }  \label{A2} \\
d\theta _{\beta }\text{ }^{\gamma } &=&\theta _{\beta }\text{ }^{\sigma
}\wedge \theta _{\sigma }\text{ }^{\gamma }  \notag
\end{eqnarray}

\noindent (Einstein summation convention used hereafter) on $H_{n},$ in
which we have used that torsion and curvature vanish on $H_{n}.$
Substituting $\theta ^{\beta }$ $=$ $\omega ^{\beta }+i\omega ^{n+\beta },$ $%
\theta _{\beta }$ $^{\gamma }$ $=$ $\omega _{\beta }$ $^{\gamma }$ $+$ $%
i\omega _{\beta }$ $^{n+\gamma }$ into (\ref{A1}) and (\ref{A2}) we obtain
the real version of structure equations: (write $\Theta $ as $\omega
^{2n+1}) $%
\begin{eqnarray}
d\omega ^{2n+1} &=&2\sum_{\beta =1}^{n}\omega ^{\beta }\wedge \omega
^{n+\beta }  \label{A3} \\
d\omega ^{\beta } &=&\omega ^{\gamma }\wedge \omega _{\gamma }\text{ }%
^{\beta }+\omega ^{n+\gamma }\wedge \omega _{n+\gamma }\text{ }^{\beta } 
\notag \\
d\omega ^{n+\beta } &=&\omega ^{\gamma }\wedge \omega _{\gamma }\text{ }%
^{n+\beta }+\omega ^{n+\gamma }\wedge \omega _{n+\gamma }\text{ }^{n+\beta }
\notag \\
d\omega _{\beta }\text{ }^{\gamma } &=&\omega _{\beta }\text{ }^{\sigma
}\wedge \omega _{\sigma }\text{ }^{\gamma }+\omega _{\beta }\text{ }%
^{n+\sigma }\wedge \omega _{n+\sigma }\text{ }^{\gamma }  \notag \\
d\omega _{\beta }\text{ }^{n+\gamma } &=&\omega _{\beta }\text{ }^{\sigma
}\wedge \omega _{\sigma }\text{ }^{n+\gamma }+\omega _{\beta }\text{ }%
^{n+\sigma }\wedge \omega _{n+\sigma }\text{ }^{n+\gamma }  \notag
\end{eqnarray}

\noindent (summation convention used in the last four lines of (\ref{A3})).
Here we have defined $\omega _{n+\gamma }$ $^{\beta }$ $:=$ $-\omega _{\beta
}$ $^{n+\gamma }$ and $\omega _{n+\gamma }$ $^{n+\beta }$ $=$ $\omega
_{\gamma }$ $^{\beta }$ so that $\omega _{a}$ $^{b}$ $=$ $-\omega _{b}$ $%
^{a} $ for $1$ $\leq $ $a,b$ $\leq $ $2n$ and $\omega _{\beta }$ $^{n+\gamma
}$ $= $ $-\omega _{n+\beta }$ $^{\gamma }$ for $1$ $\leq $ $\alpha ,\beta $ $%
\leq $ $n$ (obtained from $\theta _{\gamma }$ $^{\beta }$ being skew
hermitian)$.$

Let $\Sigma $ be a hypersurface in $H_{n}.$ Recall $\xi ^{\prime }$ $:=$ $%
(\xi \cap T\Sigma )\cap J(\xi \cap T\Sigma )$ (see (\ref{1.1})). Take $e_{n}$
$\in $ $\xi \cap T\Sigma $ $\cap $ $(\xi ^{\prime })^{\perp }$ of unit
length with respect to the Levi metric $G$ $:=$ $\frac{1}{2}d\Theta (\cdot
,J\cdot )$ $=$ $\sum_{j=1}^{n}[(dx_{j})^{2}+(dy_{j})^{2}]$ defined on $\xi .$
Let $e_{2n}$ $:=$ $Je_{n}.$ Take an orthonormal (w.r.t. $G)$ frame $e_{j},$ $%
e_{n+j},$ $1$ $\leq $ $j$ $\leq $ $n-1$ in $\xi ^{\prime }.$ Let $\{\omega
^{1},$ $..,$ $\omega ^{2n},$ $\omega ^{2n+1}$ $=$ $\Theta \}$ be the coframe
dual to $\{e_{1},$ $..,$ $e_{2n},$ $T\}.$ Recall that the function $\alpha $
on $\Sigma $ is defined so that $\alpha e_{2n}+T$ $\in $ $T\Sigma .$ Let $%
\hat{e}_{j}$ :$=$ $e_{j},$ $\hat{e}_{n+j}$ :$=$ $e_{n+j},$ $1$ $\leq $ $j$ $%
\leq $ $n-1,$ $\hat{e}_{n}$ :$=$ $e_{n},$ and $\hat{e}_{2n}$ $:=$ $\frac{%
\alpha e_{2n}+T}{\sqrt{1+\alpha ^{2}}}$ be an orthonormal basis on $T\Sigma $
with respect to the metric induced from the left invariant metric $\Theta
^{2}+G$ of $H_{n}.$ Let $\hat{\omega}^{j},$ $1$ $\leq $ $j$ $\leq $ $2n$
denote the dual coframe. Then $\omega ^{j}$ and $\hat{\omega}^{j}$ are
related as follows:%
\begin{eqnarray*}
\omega ^{j} &=&\hat{\omega}^{j}\text{ for }1\leq j\leq 2n-1 \\
\omega ^{2n} &=&\frac{\alpha }{\sqrt{1+\alpha ^{2}}}\hat{\omega}^{2n} \\
\omega ^{2n+1}( &=&\Theta )=\frac{1}{\sqrt{1+\alpha ^{2}}}\hat{\omega}^{2n}
\end{eqnarray*}

\noindent on $T\Sigma .$ The Levi-Civita connection forms $\hat{\omega}_{a}$ 
$^{b}$ are also related to pseudohermitian connection forms $\omega _{a}$ $%
^{b}$ (see \cite{CL} for more details). Define the second fundamental form $%
II_{\xi }$ $:$ $\xi \cap T\Sigma $ $\times $ $\xi \cap T\Sigma $ $%
\rightarrow $ $R$ by%
\begin{equation*}
II_{\xi }(X,Y)=-<\nabla _{Y}e_{2n},X>
\end{equation*}

\noindent where we use $<\cdot ,\cdot >$ to denote the Levi metric $G.$
Define $h_{ab}$ for $1$ $\leq $ $a,b$ $\leq $ $2n-1$ by%
\begin{equation*}
h_{ab}:=II_{\xi }(e_{a},e_{b}).
\end{equation*}

\noindent In terms of differential forms, we can write%
\begin{eqnarray}
\omega _{a}\text{ }^{2n} &=&\sum_{b=1}^{2n-1}h_{ab}\omega ^{b}+(\frac{\hat{e}%
_{a}\alpha +2\alpha ^{2}\delta _{an}}{\sqrt{1+\alpha ^{2}}})\hat{\omega}^{2n}
\label{2-0} \\
&=&\sum_{b=1}^{2n-1}h_{ab}\omega ^{b}+(e_{a}\alpha +2\alpha ^{2}\delta
_{an})\Theta .  \notag
\end{eqnarray}

\noindent by Proposition 5.5 in \cite{CL} and $\hat{\omega}^{2n}$ $=$ $\sqrt{%
1+\alpha ^{2}}\Theta $. Here $\delta _{an}$ denotes Dirac's delta function.
It is not hard to see that $h_{nn}$ is nothing but $l$ in Section 1. Note
that $II_{\xi }$ is partially symmetric, but not symmetric in general as
shown below.

\bigskip

\textbf{Proposition 2.1}. $h_{ab}$\textit{\ }$=$\textit{\ }$h_{ba}$\textit{\
for }$1$\textit{\ }$\leq $\textit{\ }$a,b$\textit{\ }$\leq $\textit{\ }$2n-1$%
\textit{\ with }$|a-b|$\textit{\ }$\neq $\textit{\ }$n$\textit{\ and }$%
h_{\beta (n+\beta )}$\textit{\ }$-$\textit{\ }$h_{(n+\beta )\beta }$\textit{%
\ }$=$\textit{\ }$2\alpha $\textit{\ for }$1$\textit{\ }$\leq $\textit{\ }$%
\beta $\textit{\ }$\leq $\textit{\ }$n-1.$

\bigskip

\proof
Observe that $\omega ^{2n}$ $-$ $\alpha \Theta $ $=0$ on $T\Sigma .$ So
using (\ref{A3}) to expand $d(\omega ^{2n}$ $-$ $\alpha \Theta )$ $=$ $0$,
we get%
\begin{eqnarray}
\sum_{c=1}^{2n-1}\omega ^{c}\wedge \omega _{c}\text{ }^{2n} &=&d\omega
^{2n}=d(\alpha \Theta )  \label{A4} \\
&=&d\alpha \wedge \Theta +\alpha d\Theta  \notag \\
&=&d\alpha \wedge \Theta +2\alpha \sum_{\beta =1}^{n}\omega ^{\beta }\wedge
\omega ^{n+\beta }  \notag
\end{eqnarray}

\noindent Applying (\ref{A4}) to $(e_{a},e_{b}),$ we obtain%
\begin{eqnarray}
h_{ab}-h_{ba} &=&\omega _{a}\text{ }^{2n}(e_{b})-\omega _{b}\text{ }%
^{2n}(e_{a})  \label{A5} \\
&=&2\alpha \sum_{\beta =1}^{n}(\delta _{a\beta }\delta _{b(n+\beta )}-\delta
_{b\beta }\delta _{a(n+\beta )}).  \notag
\end{eqnarray}

\noindent Here $\delta _{a\beta }$ denotes Dirac's delta function. The
conclusion follows from (\ref{A5}).

\endproof%

\bigskip

The results in Proposition 2.1 also appeared in \cite{CL} where a different
proof was given. Define $J^{\prime }$ on $\xi \cap T\Sigma $ by 
\begin{equation}
J^{\prime }=J\text{ on }\xi ^{\prime }\text{ and }J^{\prime }e_{n}=0.
\label{A5.1}
\end{equation}%
\noindent We can now define a shape operator $\mathfrak{S}$ $\mathfrak{:}$ $%
\xi \cap T\Sigma $ $\rightarrow $ $\xi \cap T\Sigma $ by%
\begin{equation}
\mathfrak{S}(v)=-\nabla _{v}e_{2n}+\alpha J^{\prime }v.  \label{A5.2}
\end{equation}

\bigskip

\textbf{Proposition 2.2}. (see also \cite{Rit}) $\mathfrak{S}$\textit{\ is
symmetric or self adjoint. I.e., }$<\mathfrak{S}(v_{1}),v_{2}>$\textit{\ }$=$%
\textit{\ }$<v_{1},\mathfrak{S}(v_{2})>$\textit{\ for }$v_{1},$\textit{\ }$%
v_{2}$\textit{\ }$\in $\textit{\ }$\xi \cap T\Sigma ,$\textit{\ where }$%
<\cdot ,\cdot >$\textit{\ denotes the Levi metric }$G.$

\bigskip

\proof
It suffices to show that%
\begin{eqnarray}
&<&-\nabla _{e_{a}}e_{2n}+\alpha J^{\prime }e_{a},e_{b}>  \label{A6} \\
&=&<e_{a,}-\nabla _{e_{b}}e_{2n}+\alpha J^{\prime }e_{b}>  \notag
\end{eqnarray}

\noindent for $1$ $\leq $ $a,$ $b$ $\leq $ $2n-1.$ Rewrite (\ref{A6}) as%
\begin{equation}
h_{ab}-h_{ba}=\alpha \{\delta _{(n+a)b}-\delta _{a(n+b)}\}  \label{A7}
\end{equation}

\noindent in which $n+a$ and $n+b$ are interpreted as integers from $0$ to $%
2n-1$ modulo $2n$ and $\delta _{(n+a)b}$ ($\delta _{a(n+b)},$ resp.) will
change sign if $n+a$ ($n+b,$ resp.) is larger than $2n.$ For instance, $%
\delta _{(2n+1)b}$ $=$ $-\delta _{1b}.$ Now observe that (\ref{A7}) is
equivalent to Proposition 2.1.

\endproof%

\bigskip

Recall that in Section 1 we define $X_{n}$ $\in $ $\xi ^{\prime }$ by%
\begin{equation*}
X_{n}:=\nabla _{e_{n}}e_{2n}+le_{n}
\end{equation*}

\noindent (cf. (\ref{1.2})). Observe that $\xi \cap T\Sigma $ $=$ $\xi
^{\prime }$ $\oplus $ $R$ $e_{n}.$

\bigskip

\textbf{Proposition 2.3}. \textit{At a regular point, }$X_{n}$\textit{\ }$=$%
\textit{\ }$0$\textit{\ if and only if }$\mathfrak{S}(\xi ^{\prime })$%
\textit{\ }$\subset $\textit{\ }$\xi ^{\prime }.$

\bigskip

\proof
For $v$ $\in $ $\xi ^{\prime },$ we compute%
\begin{eqnarray}
&<&\mathfrak{S}(v),e_{n}>=<-\nabla _{v}e_{2n}+\alpha Jv,e_{n}>  \label{2.8}
\\
&=&<-\nabla _{v}e_{2n},e_{n}>\text{ \ (since }Jv\in \xi ^{\prime })  \notag
\\
&=&<e_{2n},\nabla _{v}e_{n}>  \notag \\
&=&<e_{2n},\nabla _{e_{n}}v+[v,e_{n}]+Tor(v,e_{n})>  \notag
\end{eqnarray}

\noindent where $Tor(v,e_{n})$ $=$ $d\Theta (v,e_{n})T$ $=$ $0.$ Note that $%
[v,e_{n}]$ $\in $ $\xi $ $\cap $ $T\Sigma .$ So $<e_{2n},$ $[v,e_{n}]>$ $=$ $%
0.$ We therefore have%
\begin{eqnarray}
&<&\mathfrak{S}(v),e_{n}>=<e_{2n},\nabla _{e_{n}}v>  \label{2.9} \\
&=&<-\nabla _{e_{n}}e_{2n},v>  \notag \\
&=&<-X_{n}+le_{n},v>  \notag \\
&=&-<X_{n},v>  \notag
\end{eqnarray}

\noindent from (\ref{2.8}). The conclusion follows from (\ref{2.9}).

\endproof%

\bigskip

\textbf{Proposition 2.4}. \textit{At an umbilic point, choose an orthonormal
basis of }$\xi \cap T\Sigma ,$\textit{\ which are also eigenvectors of }$%
\mathfrak{S}$\textit{\ as in Proposition 1.1. Then }$h_{jm}$ $=$ $0$ \textit{%
for }$1$\textit{\ }$\leq $\textit{\ }$j,m$\textit{\ }$\leq $\textit{\ }$2n-1$%
\textit{\ except }$j$\textit{\ }$=$\textit{\ }$m$\textit{\ and }$|j-m|$%
\textit{\ }$=$\textit{\ }$n.$\textit{\ Moreover, }$h_{jj}$\textit{\ }$=$%
\textit{\ }$k$\textit{\ for }$1$\textit{\ }$\leq $\textit{\ }$j$\textit{\ }$%
\leq $\textit{\ }$2n-1,$\textit{\ }$j$\textit{\ }$\neq $\textit{\ }$n,$%
\textit{\ }$h_{nn}$\textit{\ }$=$\textit{\ }$l,$\textit{\ and }$h_{\beta
(n+\beta )}$\textit{\ }$=$\textit{\ }$-h_{(n+\beta )\beta }$\textit{\ }$=$%
\textit{\ }$\alpha $\textit{\ for }$1$\textit{\ }$\leq $\textit{\ }$\beta $%
\textit{\ }$\leq $\textit{\ }$n-1.$ \textit{In summary we can write }$h_{ab}$%
\textit{\ }$=$\textit{\ }$\omega _{a}$\textit{\ }$^{2n}(e_{b})$\textit{\ }$=$%
\textit{\ }$k\delta _{ab}+\alpha \delta _{n+a,b},$\textit{\ }$a\neq n,$%
\textit{\ }$1\leq b\leq 2n-1$\textit{\ and }$h_{nn}$\textit{\ }$=$\textit{\ }%
$l.$

\textit{\bigskip }

\proof
Compute%
\begin{equation}
<\mathfrak{S}(e_{n+\beta }),e_{\beta }>=\lambda _{n+\beta }<e_{n+\beta
},e_{\beta }>=0.  \label{2.10}
\end{equation}

\noindent On the other hand, $\mathfrak{S}(e_{n+\beta })$ $=$ $-\nabla
_{e_{n+\beta }}e_{2n}+\alpha Je_{n+\beta }$ $=$ $-\nabla _{e_{n+\beta
}}e_{2n}-\alpha e_{\beta },$ and hence%
\begin{equation}
<\mathfrak{S}(e_{n+\beta }),e_{\beta }>=h_{\beta (n+\beta )}-\alpha .
\label{2.11}
\end{equation}

\noindent The conclusion follows from (\ref{2.10}), (\ref{2.11}), and
Proposition 2.1.

\endproof%

\bigskip

The $p($or horizontal)-mean curvature $H$ of $\Sigma $ at a regular point is
defined by%
\begin{equation*}
H=\sum_{a=1}^{2n-1}h_{aa}.
\end{equation*}

\noindent Suppose $\Sigma $ is the boundary of a domain $\Omega $ in $H_{n}.$
We usually take $e_{n}$ such that the horizontal normal $e_{2n}$ $=$ $Je_{n}$
points inwards to $\Omega .$ The resulting $p$-mean curvature for a Pansu
sphere is then positive (see Example 3.2). At an umbilic point, we have%
\begin{equation*}
H=l+(2n-2)k
\end{equation*}

\noindent by Proposition 2.4.

\bigskip

\section{Umbilicity and examples}

\textbf{Proposition 3.1. }\textit{If }$\Sigma $\textit{\ is rotationally
symmetric, then it is umbilic. If, in addition, it is closed and satisfies
the condition }$l=2k$\textit{\ , then }$\Sigma $\textit{\ must be the Pansu
sphere }$S_{\lambda }$ with $\lambda $ $=$ $k$\textit{.}

\bigskip

\proof
Since $\Sigma $ is rotationally symmetric, it can be defined by the union of
the graphs of functions $f,-f$, where $f>0$ only depends on $|z|$ $:=$ $%
(\sum_{\beta =1}^{n}(x_{\beta }^{2}+y_{\beta }^{2}))^{1/2}$ and is defined
on a close interval $|z|\leq \rho $ for some positive constant $\rho $.
Write $t^{2}=f(|z|^{2})$. Then $u=f(|z|^{2})-t^{2}$ is a defining function.
We choose $e_{2n}$ $:=$ $\frac{\nabla _{b}u}{|\nabla _{b}u|}$ as the
horizontal normal so that $e_{n}$ $:=-Je_{2n}$ defines the one-dimensional
foliation on the regular part of $\Sigma $. On the regular part, we have 
\begin{equation}
\begin{split}
e_{2n}& =\sum_{\beta =1}^{n}\frac{(f^{\prime }x_{\beta }-ty_{\beta })%
\mathring{e}_{\beta }+(f^{\prime }y_{\beta }+tx_{\beta })\mathring{e}%
_{n+\beta }}{|z|\sqrt{(f^{\prime })^{2}+f}} \\
e_{n}& =\sum_{\beta =1}^{n}\frac{(f^{\prime }y_{\beta }+tx_{\beta })%
\mathring{e}_{\beta }-(f^{\prime }x_{\beta }-ty_{\beta })\mathring{e}%
_{n+\beta }}{|z|\sqrt{(f^{\prime })^{2}+f}},
\end{split}
\label{fore2n}
\end{equation}%
where $f=f(r),f^{\prime }=f^{\prime }(r)$ and $r=|z|^{2}$. Since $|\nabla
_{b}u|=2|z|\sqrt{(f^{\prime })^{2}+f}$, we see that the north pole and south
pole are the only singular points of $\Sigma $, that is, those points at $%
|z|=0$.

In order to prove that $\Sigma $ is umbilic, we are going to compute the
covariant derivatives $\nabla _{e_{n}}e_{2n}$ and $\nabla _{e}e_{2n}$, for
all $e\in \xi ^{\prime }$. By rotational symmetry, it suffices to do the
computation at such a point $p=(z,t)=(x_{1},0,\cdots ,0,t)$, i.e., $%
z_{1}=x_{1},y_{1}=0,z_{\beta }=0,\ \text{for all}\ 2\leq \beta \leq n$. We
also assume $x_{1}>0$. Let $e=\sum_{\beta =1}^{n}(a^{\beta }\mathring{e}%
_{\beta }+a^{n+\beta }\mathring{e}_{n+\beta })$, then 
\begin{equation}
\begin{split}
& \ e\in \xi ^{\prime }(p) \\
\Leftrightarrow & \ e\perp e_{2n}\ \text{and}\ e\perp e_{n} \\
\Leftrightarrow & 
\begin{array}{cc}
a^{1}f^{\prime }x_{1}+a^{n+1}tx_{1} & =0 \\ 
a^{1}tx_{1}-a^{n+1}f^{\prime }x_{1} & =0%
\end{array}
\\
\Leftrightarrow & \ a^{1}=a^{n+1}=0 \\
\Leftrightarrow & \ e=\sum_{\beta =2}^{n}(a^{\beta }\mathring{e}_{\beta
}+a^{n+\beta }\mathring{e}_{n+\beta }).
\end{split}%
\end{equation}%
Thus, if we let $e_{\beta }=\mathring{e}_{\beta +1}(p),\ e_{n+\beta }=%
\mathring{e}_{n+\beta +1}(p)$, then $\{e_{\beta },e_{n+\beta }|\ 1\leq \beta
\leq n-1\}$ constitutes an orthonormal basis of $\xi ^{\prime }(p)$. From
the formula (\ref{fore2n}) for the horizontal normal $e_{2n}$, and note that 
$\nabla \mathring{e}_{\beta }=0$, we have, replacing $x_{1}$ with $|z|$, 
\begin{equation}
\begin{split}
-\nabla _{e_{\beta }}e_{2n}& =\frac{-f^{\prime }}{|z|\sqrt{(f^{\prime
})^{2}+f}}e_{\beta }-\frac{t}{|z|\sqrt{(f^{\prime })^{2}+f}}e_{n+\beta } \\
-\nabla _{e_{n+\beta }}e_{2n}& =\frac{t}{|z|\sqrt{(f^{\prime })^{2}+f}}%
e_{\beta }+\frac{-f^{\prime }}{|z|\sqrt{(f^{\prime })^{2}+f}}e_{n+\beta } \\
-\nabla _{e_{n}}e_{2n}& =\left( \frac{(|z|^{2}-f^{\prime })}{|z|\sqrt{%
(f^{\prime })^{2}+f}}-\frac{(1+2f^{\prime \prime })f|z|}{\left( (f^{\prime
})^{2}+f\right) ^{\frac{3}{2}}}\right) e_{n},\ \ \text{i.e.}\ X_{n}=0.
\end{split}
\label{3.3}
\end{equation}%
In particular, we have 
\begin{equation}
h_{\beta (n+\beta )}=\frac{t}{|z|\sqrt{(f^{\prime })^{2}+f}},\ \ \text{and}\
\ h_{(n+\beta )\beta }=-\frac{t}{|z|\sqrt{(f^{\prime })^{2}+f}}.
\end{equation}%
On the other hand, by Proposition 2.1, we have $h_{\beta (n+\beta )}-\alpha
=h_{(n+\beta )\beta }+\alpha $. It follows that 
\begin{equation}
\alpha =\frac{t}{|z|\sqrt{(f^{\prime })^{2}+f}},  \label{3.4.1}
\end{equation}
\noindent and hence 
\begin{equation}
k=\frac{-f^{\prime }}{|z|\sqrt{(f^{\prime })^{2}+f}}.
\end{equation}%
So we have shown that "rotationally symmetric" implies "umbilic". Now
suppose $l=2k$. Then from the second equation of (\ref{parintcon}), we have $%
e_{n}k=0$. Note that $e_{n}$ is never generated by the distribution $\xi
^{\prime }$ (see Proposition 4.3). Hence $k$ is a constant, say $k=\lambda $%
. We would like to solve the ODE 
\begin{equation}
\frac{-f^{\prime }}{|z|\sqrt{(f^{\prime })^{2}+f}}=\lambda .  \label{ode1}
\end{equation}%
Taking the square of both sides of (\ref{ode1}), we have 
\begin{equation}
(f^{\prime })^{2}=\lambda ^{2}r\left( (f^{\prime })^{2}+f\right) ,
\end{equation}%
hence 
\begin{equation}
(f^{\prime })^{2}=\frac{\lambda ^{2}rf}{1-\lambda ^{2}r}
\end{equation}

\noindent It follows that 
\begin{equation}
f^{\prime }=-\sqrt{\frac{\lambda ^{2}rf}{1-\lambda ^{2}r}},\ \ \ \text{for}\
\ r\leq \frac{1}{\lambda ^{2}}.  \label{ode2}
\end{equation}%
Write (\ref{ode2}) as%
\begin{equation}
\frac{df}{\sqrt{f}}=-\sqrt{\frac{\lambda ^{2}r}{1-\lambda ^{2}r}}dr,
\end{equation}%
Integrating gives 
\begin{equation}
t=f^{\frac{1}{2}}=\frac{1}{2\lambda ^{2}}\left( \lambda |z|\sqrt{1-\lambda
^{2}|z|^{2}}+\cos ^{-1}{(\lambda |z|)}\right) +C,\ \ \ |z|\leq \frac{1}{%
\lambda }
\end{equation}%
Since $0=f(\frac{1}{\lambda })$, we have $C=0$. We have shown that $\Sigma $
is the Pansu sphere $S_{\lambda }$.

\bigskip 
\endproof%

\textbf{Example 3.2}. Recall that for any $\lambda >0$, the Pansu sphere $%
S_{\lambda }$ is the union of the graphs of the functions $f$ and $-f$,
where 
\begin{equation}
f(z)=\frac{1}{2\lambda ^{2}}\left( \lambda |z|\sqrt{1-\lambda ^{2}|z|^{2}}%
+\cos ^{-1}{\lambda |z|}\right) ,\ \ \ |z|\leq \frac{1}{\lambda }.
\end{equation}%
\noindent We take the defining function $u=f(z)-t$, and $e_{2n}=\frac{\nabla
_{b}u}{|\nabla _{b}u|},\ e_{n}=-Je_{2n}$, and $e_{1},\cdots
,e_{n-1},e_{n+1},\cdots ,e_{2n-1}$ is any orthonormal frame of $\xi ^{\prime
}$. Then by (\ref{3.3}) we have, for $\beta =1,\cdots ,n-1$, 
\begin{equation}
\begin{split}
-\nabla _{e_{\beta }}e_{2n}& =\lambda e_{\beta }-\frac{\sqrt{1-\lambda
^{2}|z|^{2}}}{|z|}e_{n+\beta } \\
-\nabla _{e_{n+\beta }}e_{2n}& =\frac{\sqrt{1-\lambda ^{2}|z|^{2}}}{|z|}%
e_{\beta }+\lambda e_{n+\beta } \\
-\nabla _{e_{n}}e_{2n}& =2\lambda e_{n},\ \ \text{i.e.}\ X_{n}=0.
\end{split}
\label{secform3}
\end{equation}%
Since $\alpha =\frac{\sqrt{1-\lambda ^{2}|z|^{2}}}{|z|}$ by (\ref{3.4.1}),
the formula (\ref{secform3}) is equivalent to 
\begin{equation}
\begin{split}
-\nabla _{e_{\beta }}e_{2n}+\alpha J^{\prime }e_{\beta }& =\lambda e_{\beta }
\\
-\nabla _{e_{n+\beta }}e_{2n}+\alpha J^{\prime }e_{n+\beta }& =\lambda
e_{n+\beta } \\
-\nabla _{e_{n}}e_{2n}& =2\lambda e_{n},
\end{split}
\label{secform4}
\end{equation}%
That is, $e_{\beta },e_{n+\beta },\beta =1,\cdots ,n-1$ are all eigenvectors
of the endomorphism $-\nabla e_{2n}+\alpha J^{\prime }$. The Pansu sphere $%
S_{\lambda }$ is hence umbilic with constant principal curvature $k=\lambda $
and constant partially normal $p$-mean curvature $l=2\lambda $. Therefore
the $p$-mean curvature $H$ $=$ $l$ $+$ $(2n-2)k$ $=$ $2n\lambda .$ Actually,
the characteristic curves in $S_{\lambda }$ are the geodesics of curvature $%
\lambda $ joining the poles.

\bigskip

\textbf{Example 3.3}. The Heisenberg sphere with radius $\rho $ is the set 
\begin{equation}
S(\rho )=\{(z,t)\in H_{n}:\ |z|^{4}+4t^{2}=\rho ^{4}\},
\end{equation}%
hence $u=\rho ^{4}-|z|^{4}-4t^{2}$ is a defining function. Choose $e_{2n}=%
\frac{\nabla _{b}u}{|\nabla _{b}u|},\ e_{n}=-Je_{2n}$, and $e_{1},\cdots
,e_{n-1},e_{n+1},\cdots ,e_{2n-1}$ being any orthonormal frame of $\xi
^{\prime }$. Then by (\ref{3.3}) we have, for $\beta =1,\cdots ,n-1$, 
\begin{equation}
\begin{split}
-\nabla _{e_{\beta }}e_{2n}& =\frac{|z|}{\rho ^{2}}e_{\beta }-\frac{2t}{\rho
^{2}|z|}e_{n+\beta } \\
-\nabla _{e_{n+\beta }}e_{2n}& =\frac{2t}{\rho ^{2}|z|}e_{\beta }+\frac{|z|}{%
\rho ^{2}}e_{n+\beta } \\
-\nabla _{e_{n}}e_{2n}& =\frac{3|z|}{\rho ^{2}}e_{n},\ \ \text{i.e.}\
X_{n}=0.
\end{split}
\label{secform1}
\end{equation}%
Since $\alpha =\frac{2t}{\rho ^{2}|z|}$ by (\ref{3.4.1}), the formula (\ref%
{secform1}) is equivalent to 
\begin{equation}
\begin{split}
-\nabla _{e_{\beta }}e_{2n}+\alpha J^{\prime }e_{\beta }& =\frac{|z|}{\rho
^{2}}e_{\beta } \\
-\nabla _{e_{n+\beta }}e_{2n}+\alpha J^{\prime }e_{n+\beta }& =\frac{|z|}{%
\rho ^{2}}e_{n+\beta } \\
-\nabla _{e_{n}}e_{2n}& =\frac{3|z|}{\rho ^{2}}e_{n}.
\end{split}
\label{secform2}
\end{equation}%
That is, $e_{\beta },e_{n+\beta },\beta =1,\cdots ,n-1$ are all eigenvectors
of the endomorphism $-\nabla e_{2n}+\alpha J^{\prime }$. We see from (\ref%
{secform2}) that the Heisenberg sphere is umbilic with $l=3k$, which is not
a constant.

Now we introduce some umbilic hypersurfaces with $\alpha =0$.

\bigskip

\textbf{Example 3.4.} Let $\Sigma ^{\ast }\subset R^{2n}$ be a hypersurface
of $R^{2n}$ which defined by $f(x,y)=0$, where $x=(x_{1}\cdots ,x_{n}),\
y=(y_{1}\cdots ,y_{n})$ and the gradient $\nabla f\neq 0$ on $\Sigma ^{\ast
} $. We define the hypersurface $\Sigma _{\Sigma ^{\ast }}$ of $H_{n}$ by 
\begin{equation}
\Sigma _{\Sigma ^{\ast }}=\Sigma ^{\ast }\times R.
\end{equation}%
Then the function $u(x,y,t)=f(x,y)$ is a defining function of $\Sigma
_{\Sigma ^{\ast }}$. We have 
\begin{equation}
e_{2n}=\frac{\sum_{\beta =1}^{n}\left( f_{\beta }\mathring{e}_{\beta
}+f_{n+\beta }\mathring{e}_{n+\beta }\right) }{\sqrt{\sum_{\beta
=1}^{n}(f_{\beta })^{2}+(f_{n+\beta })^{2}}},
\end{equation}%
where $f_{\beta }=\frac{\partial f}{\partial x_{\beta }},\ f_{n+\beta }=%
\frac{\partial f}{\partial y_{\beta }}$. Since both $T=\frac{\partial }{%
\partial t}$ and $\alpha e_{2n}+T$ are tangent to $\Sigma _{\Sigma ^{\ast }}$%
, we see that $\alpha =0$ on $\Sigma _{\Sigma ^{\ast }}$.

(1) Suppose $f(x,y)=\sum_{\beta =1}^{n}A^{\beta }x_{\beta }+A^{n+\beta
}y_{\beta }$. Then $f(x,y)=0$ defines a hyperplane $E$ in $R^{2n}$. We have 
\begin{equation}
e_{2n}=\frac{\sum_{a=1}^{2n}A^{a}\mathring{e}_{a}}{\sqrt{%
\sum_{a=1}^{2n}(A^{a})^{2}}}.
\end{equation}%
Therefore we have $\nabla e_{2n}=0$. Since $\alpha =0$, this implies that
the hypersurface $\Sigma _{E}$ in $H_{n}$ is umbilic with $l=k=0$.

(2) Suppose $f(x,y)=\sum_{\beta =1}^{n}\left( x_{\beta }^{2}+y_{\beta
}^{2}\right) -c^{2}$ for some constant $c>0$. Then $f(x,y)=0$ defines a $%
(2n-1)$-dimensional sphere $S^{2n-1}(c)$ in $R^{2n}$ with radius $c$. If we
choose $u(x,y,t)=-f(x,y)$, then 
\begin{equation}
e_{2n}=-\frac{\sum_{\beta =1}^{n}\left( x_{\beta }\mathring{e}_{\beta
}+y_{\beta }\mathring{e}_{n+\beta }\right) }{\sqrt{\sum_{\beta
=1}^{n}(x_{\beta })^{2}+(y_{\beta })^{2}}}=-\frac{\sum_{\beta =1}^{n}\left(
x_{\beta }\mathring{e}_{\beta }+y_{\beta }\mathring{e}_{n+\beta }\right) }{c}%
.
\end{equation}%
For any $X=\sum_{\beta =1}^{n}\left( a^{\beta }\mathring{e}_{\beta
}+a^{n+\beta }\mathring{e}_{n+\beta }\right) \in T\Sigma \cap \xi $, we have 
\begin{equation}
\begin{split}
-\nabla _{X}e_{2n}& =\frac{1}{c}\nabla \left( \sum_{\beta =1}^{n}(Xx_{\beta
})\mathring{e}_{\beta }+(Xy_{\beta })\mathring{e}_{n+\beta }\right) \\
& =\frac{1}{c}\sum_{\beta =1}^{n}\left( a^{\beta }\mathring{e}_{\beta
}+a^{n+\beta }\mathring{e}_{n+\beta }\right) \\
& =\frac{1}{c}X.
\end{split}%
\end{equation}%
Since $\alpha =0$, this implies that the hypersurface $\Sigma _{S^{2n-1}(c)}$
is umbilic with $l=k=\frac{1}{c}$.

\section{Properties of umbilic hypersurfaces}

\textbf{Proposition 4.1}. \textit{Suppose }$\Sigma $\textit{\ is an umbilic
hypersurface. If }$p\in \Sigma $\textit{\ is a singular point, then it is
isolated.}

\bigskip

\proof
After the action of the left translation $L_{p^{-1}}$, locally around $p$,
the hypersurface can be represented by the graph of a function $t=u(x,y)$
defined on a domain $\Omega \subset R^{2n}$ with $(0,0)\in \Omega ,\
u(0,0)=0,\ u_{x_{\beta }}(0,0)=u_{y_{\beta }}(0,0)=0$, where $%
x=(x_{1},\cdots ,x_{n}),\ y=(y_{1}\cdots ,y_{n})$. Moreover, after a
suitable orthogonal transformation on $R^{2n}$, we can assume, without loss
of generality, that the function $z=u(x,y)$ has the canonical diagonal forms 
\begin{equation}
t=u(x,y)=\sum_{\beta =1}^{n}\left( B_{\beta }x_{\beta }^{2}+B_{n+\beta
}x_{n+\beta }^{2}\right) +O(3),  \label{canoform}
\end{equation}%
for some constants $B_{\beta },B_{n+\beta }$, where we sometimes use $%
x_{n+\beta }$ instead of $y_{\beta }$. Consider the map $\varphi $ $:$ $q$ $%
\in $ $\Omega $ $\rightarrow $ $(\nabla u+\vec{F})(q)$ $\in $ $R^{2n}$ where 
$\vec{F}$ $:=$ $(-y_{1},$ $..,$ $-y_{n},$ $x_{1},$ $..,$ $x_{n}).$ To show
that $p$ $(=(0,0,0))$ is isolated, it is sufficient to show $\ker d\varphi
((0,0))$ $=$ $\{0\}$ by the implicit function theorem. So in matrix form, it
is sufficient to show that the following $(2n\times 2n)$-matrix is of full
rank 
\begin{equation}
U(p)=\left[ 
\begin{array}{cc}
u_{\beta \gamma } & u_{\beta (n+\gamma )} \\ 
u_{(n+\beta )\gamma } & u_{(n+\beta )(n+\gamma )}%
\end{array}%
\right] (0,0)+\left[ 
\begin{array}{cc}
0 & -I_{n} \\ 
I_{n} & 0%
\end{array}%
\right] .  \label{kma}
\end{equation}%
It is easy to see that 
\begin{equation}
\begin{split}
u_{b}& =\frac{\partial u}{\partial x_{b}}=2B_{b}x_{b}+O(2), \\
u_{ba}& =2B_{b}\delta _{ba}+O(1),\ \ \text{for}\ 1\leq a,b\leq 2n.
\end{split}%
\end{equation}%
Hence we have 
\begin{equation}
U(p)=\left[ 
\begin{array}{cccccc}
2B_{1} & \cdots & 0 &  &  &  \\ 
\vdots & \ddots & \vdots &  & -I_{n} &  \\ 
0 & \cdots & 2B_{n} &  &  &  \\ 
&  &  & 2B_{n+1} & \cdots & 0 \\ 
& I_{n} &  & \vdots & \ddots & \vdots \\ 
&  &  & 0 & \cdots & 2B_{2n}%
\end{array}%
\right] .
\end{equation}%
We will show that if $\Sigma $ is umbilic, then $B_{1}=B_{2}=\cdots =B_{2n}$
(write this common value as $B$). So it follows from basic linear algebra
that the determinant of $U(p)$ equals $(4B^{2}+1)^{n}$ $\neq $ $0$. The
matrix $U(p)$ is therefore of full rank (another argument is to observe that
the kernel of $U(p)$ as a linear transformation consists of zero vector
only), which implies that $p$ is isolated. Let 
\begin{equation*}
\rho =u(x,y)-t,
\end{equation*}%
which is a defining function. We have 
\begin{equation}
\begin{split}
\nabla _{b}\rho & =\sum_{a=1}^{2n}(\mathring{e}_{a}\rho )\mathring{e}_{a} \\
& =\sum_{\beta =1}^{n}(u_{\beta }-y_{\beta })\mathring{e}_{\beta
}+(u_{n+\beta }+x_{\beta })\mathring{e}_{n+\beta },
\end{split}%
\end{equation}%
and hence 
\begin{equation}
e_{2n}=\frac{\nabla _{b}\rho }{|\nabla _{b}\rho |}=\frac{\sum_{\beta
=1}^{n}(u_{\beta }-y_{\beta })\mathring{e}_{\beta }+(u_{n+\beta }+x_{\beta })%
\mathring{e}_{n+\beta }}{D},
\end{equation}%
where 
\begin{equation}
D=\sqrt{(u_{\beta }-y_{\beta })^{2}+(u_{n+\beta }+x_{\beta })^{2}}.
\end{equation}%
Since $\Sigma $ is umbilic, for any $e=a^{\beta }\mathring{e}_{\beta
}+a^{n+\beta }\mathring{e}_{n+\beta }\in \xi ^{\prime }$, we have 
\begin{equation}
-\nabla _{e}e_{2n}+\alpha J^{\prime }e=ke,  \label{umeq1}
\end{equation}%
where $k$ is the common eigenvalue of the operator $-\nabla e_{2n}+\alpha
J^{\prime }$. From (\ref{umeq1}). For any $e\in \xi ^{\prime },\ |e|=1$, we
compute 
\begin{equation}
\begin{split}
k& =\big<ke,e\big>=\big<-\nabla _{e}e_{2n},e\big> \\
& =-a^{\beta }e\left( \frac{u_{\beta }-y_{\beta }}{D}\right) -a^{n+\beta
}e\left( \frac{u_{n+\beta }+x_{\beta }}{D}\right) \\
& =-a^{\beta }\frac{e(u_{\beta }-y_{\beta })}{D}-a^{n+\beta }\frac{%
e(u_{n+\beta }+x_{\beta })}{D} \\
& \ \ \ \ -\left( a^{\beta }(u_{\beta }-y_{\beta })+a^{n+\beta }(u_{n+\beta
}+x_{\beta })\right) e\left( \frac{1}{D}\right) \\
& =-a^{\beta }\frac{e(u_{\beta }-y_{\beta })}{D}-a^{n+\beta }\frac{%
e(u_{n+\beta }+x_{\beta })}{D},
\end{split}
\label{umeq2}
\end{equation}%
where for the last equality, we have used the fact that $e\in \xi ^{\prime }$%
, and hence 
\begin{eqnarray*}
0 &=&\big<e,e_{2n}\big> \\
&=&\frac{1}{D}\left( a^{\beta }(u_{\beta }-y_{\beta })+a^{n+\beta
}(u_{n+\beta }+x_{\beta })\right) .
\end{eqnarray*}%
Now we compute 
\begin{equation}
\begin{split}
e(u_{\beta }-y_{\beta })& =a^{\gamma }\mathring{e}_{\gamma }(u_{\beta
}-y_{\beta })+a^{n+\gamma }\mathring{e}_{n+\gamma }(u_{\beta }-y_{\beta }) \\
& =2B_{\beta }a^{\beta }-a^{n+\beta }+O(1),
\end{split}
\label{umeq3}
\end{equation}%
and 
\begin{equation}
\begin{split}
e(u_{n+\beta }+x_{\beta })& =a^{\gamma }\mathring{e}_{\gamma }(u_{n+\beta
}+x_{\beta })+a^{n+\gamma }\mathring{e}_{n+\gamma }(u_{n+\beta }+x_{\beta })
\\
& =2B_{n+\beta }a^{n+\beta }+a^{\beta }+O(1).
\end{split}
\label{umeq4}
\end{equation}%
where O(1) means a function bounded by constant times $r^{1}$ ($%
=(\tsum\limits_{\beta =1}^{n}(x_{\beta }^{2}+y_{\beta }^{2}))^{1/2})$ when
evaluate in a small neighborhood of the origin. Substituting (\ref{umeq3})
and (\ref{umeq4}) into (\ref{umeq2}), we get, for any fixed regular point $q$
(in a small neighborhood of the origin), 
\begin{equation}
k=\frac{-2(a^{\beta })^{2}B_{\beta }-2(a^{n+\beta })^{2}B_{n+\beta }+O(1)}{D}%
,
\end{equation}%
or 
\begin{equation}
2\Big(a^{\beta }(q)\Big)^{2}B_{\beta }+2\Big(a^{n+\beta }(q)\Big)%
^{2}B_{n+\beta }=\Big(-kD+O(1)\Big)(q),  \label{umeq5}
\end{equation}%
for any $e\in \xi ^{\prime }$ with $|e|=1$. Since $\sum_{\beta =1}^{n}\Big(%
a^{\beta }(q)\Big)^{2}+\Big(a^{n+\beta }(q)\Big)^{2}=|e|^{2}=1$, the left
hand side of (\ref{umeq5}) is just the average value of $B_{\beta },\
B_{n+\beta },\ 1\leq \beta \leq n$, with weight $\Big(a^{\beta }(q)\Big)%
^{2},\ \Big(a^{n+\beta }(q)\Big)^{2},$ respectively. On the other hand, we
see that the right hand side is a constant (independent of $a^{\beta },$ $%
a^{n+\beta })$ for a fixed regular point $q$. Therefore formula (\ref{umeq5}%
) means that the average value of $B_{\beta },\ B_{n+\beta },\ 1\leq \beta
\leq n$, for any weight is a constant. Notice that the space of all weights
is a sphere with dimension $2n-3$, which is positive for $n\geq 2$. This
implies that $B_{1}=B_{2}=\cdots =B_{2n}$.

\endproof%

\bigskip

\textbf{Proposition 4.2. }\textit{Suppose }$\Sigma $\textit{\ is an umbilic
hypersurface. Then we have} 
\begin{equation}
\begin{split}
ek& =el=e\alpha =e(e_{n}\alpha )=0,\ \ \text{for all}\ e\in \xi ^{\prime },
\\
e_{n}k& =(l-2k)\alpha ,\ \ \ \ \ \ \ \hat{e}_{2n}k=\frac{\alpha
(k^{2}+e_{n}\alpha +\alpha ^{2})}{\sqrt{1+\alpha ^{2}}}, \\
e_{n}\alpha & =k^{2}-\alpha ^{2}-kl,\ \ \ \hat{e}_{2n}\alpha =\frac{%
-k(e_{n}\alpha )}{\sqrt{1+\alpha ^{2}}}, \\
\hat{e}_{2n}l& =\frac{e_{n}e_{n}\alpha +6\alpha e_{n}\alpha +4\alpha
^{3}+\alpha l^{2}}{\sqrt{1+\alpha ^{2}}}
\end{split}
\label{parintcon}
\end{equation}

\bigskip

The proof of Proposition 4.2 is a tedious computation. We will show the
computation in Section 5.

\bigskip

\textbf{Proposition 4.3}. \textit{Suppose }$\Sigma $\textit{\ is an umbilic
hypersurface. Let }$L(\xi ^{\prime })$\textit{\ denote the smallest }$%
C^{\infty }$\textit{-module which contains }$e_{1},\cdots
,e_{n-1},e_{n+1},\cdots ,e_{2n-1}$\textit{\ and is closed under the Lie
bracket. Then the rank of }$L(\xi ^{\prime })$\textit{\ is }$2n-1$\textit{.
Therefore, by Frobenius theorem, the module defines a }$(2n-1)$\textit{%
-dimensional foliation. Moreover, the characteristic direction }$e_{n}$%
\textit{\ is always transversal to each leaf of the }$(2n-1)$\textit{%
-dimensional foliation.}

\bigskip

\proof
For $1\leq \beta \leq n$, let 
\begin{equation}
Z_{\beta }=\frac{1}{2}(e_{\beta }-ie_{n+\beta }).
\end{equation}%
Let $Z_{\bar{\beta}}$ denote the complex conjugate of $Z_{\beta }.$ We claim 
\begin{equation}
\begin{split}
\lbrack Z_{\beta },Z_{\gamma }]& =0,\ \ \text{mod}\ \ Z_{\sigma },1\leq
\sigma \leq n-1, \\
\lbrack Z_{\bar{\beta}},Z_{\gamma }]& =0,\ \ \text{mod}\ \ Z_{\sigma },Z_{%
\bar{\sigma}},1\leq \sigma \leq n-1,\ \ \text{for}\ \ \beta \neq \gamma .
\end{split}
\label{lie1}
\end{equation}%
and 
\begin{equation}
\begin{split}
\lbrack Z_{\bar{\beta}},Z_{\beta }]& =iT+i\alpha e_{2n}+ike_{n} \\
& =i\sqrt{1+\alpha ^{2}}\hat{e}_{2n}+ike_{n},\ \ \text{mod}\ \ Z_{\sigma
},Z_{\bar{\sigma}},1\leq \sigma \leq n-1.
\end{split}
\label{lie2}
\end{equation}%
Finally, for each $\beta $, with $1\leq \beta \leq n-1$, we also claim 
\begin{equation}
\lbrack Z_{\beta },(iT+i\alpha e_{2n}+ike_{n})]=0,\ \ \text{mod}\ \
Z_{\sigma },Z_{\bar{\sigma}},1\leq \sigma \leq n-1.  \label{lie3}
\end{equation}%
From (\ref{lie1}), (\ref{lie2}) and (\ref{lie3}), we see that the rank of $%
L(\xi ^{\prime })$ is $2n-1$. In particular, from (\ref{lie2}) and (\ref%
{lie3}), we see that the distribution never generates the direction $e_{n}$.
In order to complete the proof, we now carry out the computation for (\ref%
{lie1}), (\ref{lie2}) and (\ref{lie3}). First we are going to show formulae (%
\ref{lie1}). For $1\leq \beta ,\gamma \leq n-1,\ \beta \neq \gamma $, we
have 
\begin{equation}
\lbrack Z_{\beta },Z_{\gamma }]=\left( \sum_{\rho =1}^{n-1}\theta _{\gamma
}{}^{\rho }(Z_{\beta })Z_{\rho }-\theta _{\beta }{}^{\rho }(Z_{\gamma
})Z_{\rho }\right) +\theta _{\gamma }{}^{n}(Z_{\beta })Z_{n}-\theta _{\beta
}{}^{n}(Z_{\gamma })Z_{n},
\end{equation}%
(see Section 4 in \cite{Lee}) where 
\begin{equation}
\begin{split}
\theta _{\gamma }{}^{n}(Z_{\beta })& =\omega _{n+\gamma }{}^{2n}(Z_{\beta
})+i\omega _{\gamma }{}^{2n}(Z_{\beta }) \\
& =\frac{1}{2}\left( \omega _{n+\gamma }{}^{2n}(e_{\beta })-i\omega
_{n+\gamma }{}^{2n}(e_{n+\beta })\right) +\frac{1}{2}i\left( \omega _{\gamma
}{}^{2n}(e_{\beta })-i\omega _{\gamma }{}^{2n}(e_{n+\beta })\right) \\
& =\frac{1}{2}\left( h_{(n+\gamma )\beta }-ih_{(n+\gamma )(n+\beta )}\right)
+\frac{1}{2}i\left( h_{\gamma \beta }-ih_{\gamma (n+\beta )}\right) \\
& =0,
\end{split}%
\end{equation}%
for the last equality, we have used the fact $h_{jk}=0,\ \text{for}\ 1\leq
j,k\leq 2n-1$, except $j=k$ or $|j-k|=n$ by Proposition 2.4. Similarly, we
have 
\begin{equation}
\theta _{\beta }{}^{n}(Z_{\gamma })=0.
\end{equation}%
Thus we have shown the first equation of (\ref{lie1}). The proof of the
second equation of (\ref{lie1}) is similar (note that the Levi metric $%
h_{\beta \bar{\gamma}}=\delta _{\beta \gamma }$). Next, we are going to show
(\ref{lie2}). For $1\leq \beta \leq n-1$, we have%
\begin{equation}
\begin{split}
\lbrack Z_{\bar{\beta}},Z_{\beta }]& =iT+\theta _{\beta }{}^{n}(Z_{\bar{\beta%
}})Z_{n}-\theta _{\bar{\beta}}{}^{\bar{n}}(Z_{\beta })Z_{\bar{n}},\ \ \text{%
mod}\ Z_{\rho },Z_{\bar{\rho}},\ 1\leq \rho \leq n-1 \\
& =iT+\left( \frac{1}{2}[(\omega _{n+\beta }{}^{2n}+i\omega _{\beta
}{}^{2n})(e_{\beta }+ie_{n+\beta })]Z_{n}-\ \text{conjugate}\right) \\
& =iT+\left( \left[ \frac{1}{2}(h_{(n+\beta )\beta }-h_{\beta (n+\beta )})+%
\frac{1}{2}i(h_{\beta \beta }+h_{(n+\beta )(n+\beta )})\right] Z_{n}-\ \text{%
conjugate}\right) \\
& =iT+\left( \frac{1}{2}(-\alpha +ik)(e_{n}-ie_{2n})-\ \text{conjugate}%
\right) \\
& =i(T+\alpha e_{2n})+ike_{n},\ \ \text{mod}\ Z_{\rho },Z_{\bar{\rho}},\
1\leq \rho \leq n-1.
\end{split}%
\end{equation}%
For the above computation, we have used the fact that $h_{\beta \beta
}=h_{(n+\beta )(n+\beta )}=k$ and $h_{\beta (n+\beta )}-h_{(n+\beta )\beta
}=2\alpha $ (see Proposition 2.4). We have shown (\ref{lie2}). Finally, we
will show (\ref{lie3}). For $1\leq \beta \leq n-1$, we also need the fact
that $Z_{\beta }\alpha =Z_{\beta }k=0$ by Proposition 4.2. Therefore we have%
\begin{equation}
\begin{split}
\lbrack Z_{\beta },(iT+i\alpha e_{2n}+ike_{n})]& =i[Z_{\beta },\alpha
e_{2n}]+i[Z_{\beta },T]+i[Z_{\beta },ke_{n}] \\
& =i\alpha \lbrack Z_{\beta },e_{2n}]+i[Z_{\beta },T]+ik[Z_{\beta },e_{n}] \\
& =-\alpha \lbrack Z_{\beta },Z_{n}-Z_{\bar{n}}]+ik[Z_{\beta },Z_{n}+Z_{\bar{%
n}}]+i[Z_{\beta },T] \\
& =-(\alpha -ik)[Z_{\beta },Z_{n}]+(\alpha +ik)[Z_{\beta },Z_{\bar{n}%
}]+i[Z_{\beta },T]
\end{split}
\label{lie33}
\end{equation}%
where 
\begin{equation}
\begin{split}
\lbrack Z_{\beta },Z_{n}]& =\theta _{n}{}^{n}(Z_{\beta })Z_{n}-\theta
_{\beta }{}^{n}(Z_{n})Z_{n} \\
& =\frac{1}{2}[(\omega _{2n}{}^{2n}+i\omega _{n}{}^{2n}(e_{\beta
}-ie_{n+\beta })]Z_{n}-\theta _{\beta }{}^{n}(Z_{n})Z_{n} \\
& =\frac{1}{2}i(h_{n\beta }-ih_{n(n+\beta )})Z_{n}-\theta _{\beta
}{}^{n}(Z_{n})Z_{n} \\
& =-\theta _{\beta }{}^{n}(Z_{n})Z_{n},\ \ \text{mod}\ Z_{\rho },1\leq \rho
\leq n-1.
\end{split}
\label{lie34}
\end{equation}%
Here we have used Proposition 2.4. Similarly, we have%
\begin{equation}
\lbrack Z_{\beta },Z_{\bar{n}}]=-\theta _{\beta }{}^{n}(Z_{\bar{n}})Z_{n},\
\ \text{mod}\ Z_{\rho },1\leq \rho \leq n-1.  \label{lie35}
\end{equation}%
Since the pseudohermitian torsion for $H_{n}$ is zero, we have 
\begin{equation}
\lbrack Z_{\beta },T]=-\theta _{\beta }{}^{n}(T)Z_{n},\ \ \text{mod}\
Z_{\rho },\ 1\leq \rho \leq n-1.  \label{lie36}
\end{equation}%
Substituting (\ref{lie34}), (\ref{lie35}) and (\ref{lie36}) into (\ref{lie33}%
), we obtain 
\begin{equation}
\begin{split}
\lbrack Z_{\beta },(iT+i\alpha e_{2n}+ike_{n})]& =(\alpha -ik)\theta _{\beta
}{}^{n}(Z_{n})Z_{n}-(\alpha +ik)\theta _{\beta }{}^{n}(Z_{\bar{n}%
})Z_{n}-i\theta _{\beta }{}^{n}(T)Z_{n} \\
& =\theta _{\beta }{}^{n}\Big((\alpha -ik)Z_{n}-(\alpha +ik)Z_{\bar{n}}\Big)%
Z_{n}-i\theta _{\beta }{}^{n}(T)Z_{n} \\
& =-i\theta _{\beta }{}^{n}(ke_{n}+\alpha e_{2n})Z_{n}-i\theta _{\beta
}{}^{n}(T)Z_{n} \\
& =-i\theta _{\beta }{}^{n}(ke_{n}+\alpha e_{2n}+T)Z_{n} \\
& =-\left( ik\theta _{\beta }{}^{n}(e_{n})+i\sqrt{1+\alpha ^{2}}\theta
_{\beta }{}^{n}(\hat{e}_{2n})\right) Z_{n},
\end{split}%
\end{equation}%
where 
\begin{equation}
\theta _{\beta }{}^{n}(e_{n})=(\omega _{n+\beta }{}^{2n}+i\omega _{\beta
}{}^{2n})(e_{n})=h_{(n+\beta )n}+ih_{\beta n}=0
\end{equation}%
by Proposition 2.4 and 
\begin{equation}
\theta _{\beta }{}^{n}(\hat{e}_{2n})=\omega _{n+\beta }{}^{2n}(\hat{e}%
_{2n})+i\omega _{\beta }{}^{2n}(\hat{e}_{2n})=\frac{e_{n+\beta }\alpha }{%
\sqrt{1+\alpha ^{2}}}+i\frac{e_{\beta }\alpha }{\sqrt{1+\alpha ^{2}}}=0,
\end{equation}%
by the first formula of (\ref{parintcon}). We have completed the proof.

\endproof%

\bigskip

From Proposition 4.2 and Proposition 4.3, we have

\bigskip

\textbf{Proposition 4.4}. \textit{Suppose }$\Sigma $\textit{\ is an umbilic
hypersurface. Then the common eigenvalue }$k$\textit{, the fundamental
function }$\alpha $\textit{\ and the partially p-mean curvature }$l$\textit{%
\ are all constants on each leaf of the foliation described in Proposition
4.3.}

\bigskip

\textbf{Proposition 4.5}. \textit{Suppose }$\Sigma $\textit{\ is umbilic and
satisfies the condition }$l=2k$\textit{. Then }$k$\textit{\ must be
constant, say }$k=\lambda $\textit{, and each characteristic curve is a
geodesic of curvature }$\lambda $\textit{. That is, the regular part of }$%
\Sigma $\textit{\ is foliated by geodesics of curvature }$\lambda $\textit{.}

\textit{\bigskip }

\proof
From the second equation of (\ref{parintcon}), we see that the condition $%
l=2k$ implies that $e_{n}k=0$. On the other hand, from Proposition 4.3 and
Proposition 4.4, we see that $k$ is constant on each leaf and $e_{n}$ is
transversal to each leaf. Thus $k$, and hence $l$, are constant on the whole
regular part of $\Sigma $, say $k=\lambda $. Therefore we have 
\begin{equation}
\begin{split}
-\nabla _{e_{n}}e_{2n}& =le_{n} \\
& =2ke_{n}=2\lambda e_{n}.
\end{split}%
\end{equation}%
This equation is equivalent to 
\begin{equation}
\nabla _{e_{n}}e_{n}-2\lambda Je_{n}=0,
\end{equation}%
which implies that each characteristic curve is a geodesic of curvature $%
\lambda $ (see page 52 in \cite{Rit}).

\endproof%

\bigskip

\proof
\textbf{(of Theorem 1.4) }From Proposition 4.5, we see that the regular part
of $\Sigma $ is foliated by geodesics of curvature $\lambda $. If $\lambda
>0 $, then $\Sigma ,$ containing the singular point $p$ which is isolated by
Proposition 4.1, is congruent with part of the Pansu sphere $S_{\lambda }$
by Heisenberg translating $p$ to a pole of $S_{\lambda }$. On the other
hand, if $\lambda $ $=$ $0$, then $\Sigma $ is congruent with part of a
hyperplane orthogonal to the $t$-axis by a similar reasoning.

\endproof%

\bigskip

\proof
\textbf{(of Theorem 1.5) }If $\alpha =0$, then from the second equation of (%
\ref{parintcon}), we have $e_{n}k=0$. On the other hand, from Proposition
4.3 and Proposition 4.4, we see that $k$ is constant on each leaf of the
filiation defined by the module $L(\xi ^{\prime })$ and the characteristic
direction $e_{n}$ is always transversal to each leaf. Therefore we conclude
that $k$ must be constant on $\Sigma $ (note that $\alpha =0$ implies that $%
\Sigma $ contains no singular point).

Next, from the third equation of (\ref{parintcon}), we have%
\begin{equation*}
0=e_{n}\alpha =k(k-l),
\end{equation*}%
\noindent which implies that $k=l$, provided that $k>0$. From this and
Proposition 4.5, we obtain that $\Sigma $ is foliated by geodesics whose
projections on the $xy$-space lie in Euclidean spheres. On the other hand, $%
\alpha =0$ implies that $\Sigma $ is a vertical hypersurface, that is, the
vertical vector $T=\frac{\partial }{\partial t}$ is always tangent to $%
\Sigma $ at each point. Therefore $\Sigma $ is congruent with part of the
hypersurface $\Sigma _{S^{2n-1}(c)}$ for some $c>0$. If $k=l=0$, a similar
argument shows that $\Sigma $ is foliated by straight lines. Since $\alpha $ 
$=$ $0$, we conclude that $\Sigma $ is congruent with part of the
hypersurface $\Sigma _{E}$ for some hyperplane $E$ in $R^{2n}$.

\endproof%

\bigskip

\section{Proof of Proposition 4.2}

In this section, we will prove Proposition 4.2. Observe that since $\Sigma $
is umbilic, we have, for $1\leq \beta \leq n-1$, 
\begin{equation}
\begin{split}
\omega _{2n}{}^{\beta }& =-k\omega ^{\beta }-\alpha \omega ^{n+\beta
}+\omega _{2n}{}^{\beta }(\hat{e}_{2n})\hat{\omega}^{2n} \\
\omega _{2n}{}^{n+\beta }& =\alpha \omega ^{\beta }-k\omega ^{n+\beta
}+\omega _{2n}{}^{n+\beta }(\hat{e}_{2n})\hat{\omega}^{2n} \\
\omega _{2n}{}^{n}& =-l\omega ^{n}+\omega _{2n}{}^{n}(\hat{e}_{2n})\hat{%
\omega}^{2n}.
\end{split}
\label{expofcon}
\end{equation}%
\noindent due to $h_{\beta \beta }$ $=$ $k,$ $h_{\beta (n+\beta )}$ $=$ $%
\alpha $ by Proposition 2.4. For each $k,1\leq k\leq 2n-1$, expanding the
following partial integrability conditions (see (\ref{A3})) 
\begin{equation}
d\omega _{k}{}^{2n}=\omega _{k}{}^{\gamma }\wedge \omega _{\gamma
}{}^{2n}+\omega _{k}{}^{n+\gamma }\wedge \omega _{n+\gamma }{}^{2n}
\end{equation}%
(summation convention used) and comparing the coefficients of the
corresponding terms, we will then get formulae (\ref{parintcon}). Now we
perform the computation. Let%
\begin{equation*}
\Gamma ^{j}:=\omega _{2n}{}^{j}(\hat{e}_{2n}),1\leq j\leq 2n-1.
\end{equation*}%
Taking the exterior differential of the first formula of (\ref{expofcon}),
we get 
\begin{equation}
\begin{split}
d\omega _{2n}{}^{\beta }& =d(-k\omega ^{\beta }-\alpha \omega ^{n+\beta
}+\Gamma ^{\beta }\hat{\omega}^{2n}) \\
& =-\sum_{a=1}^{2n}(\hat{e}_{a}k)\hat{\omega}^{a}\wedge \omega ^{\beta
}-kd\omega ^{\beta }-\sum_{a=1}^{2n}(\hat{e}_{a}\alpha )\hat{\omega}%
^{a}\wedge \omega ^{n+\beta }-\alpha d\omega ^{n+\beta } \\
& \ \ \ +\sum_{a=1}^{2n}(\hat{e}_{a}\Gamma ^{\beta })\hat{\omega}^{a}\wedge 
\hat{\omega}^{2n}+\Gamma ^{\beta }d\hat{\omega}^{2n} \\
& =-\sum_{a=1}^{2n}(\hat{e}_{a}k)\hat{\omega}^{a}\wedge \omega ^{\beta
}-k\sum_{a=1}^{2n}\hat{\omega}^{a}\wedge \hat{\omega}_{a}{}^{\beta
}-\sum_{a=1}^{2n}(\hat{e}_{a}\alpha )\hat{\omega}^{a}\wedge \omega ^{n+\beta
} \\
& \ \ \ -\alpha \sum_{a=1}^{2n}\hat{\omega}^{a}\wedge \hat{\omega}%
_{a}{}^{n+\beta }+\sum_{a=1}^{2n}(\hat{e}_{a}\Gamma ^{\beta })\hat{\omega}%
^{a}\wedge \hat{\omega}^{2n}+\Gamma ^{\beta }\sum_{a=1}^{2n}\hat{\omega}%
^{a}\wedge \hat{\omega}_{a}{}^{2n},
\end{split}
\label{dofcon1}
\end{equation}%
where we have used the following formulae 
\begin{equation}
\begin{split}
& d\omega ^{\beta }=d\hat{\omega}^{\beta }=\sum_{a=1}^{2n}\hat{\omega}%
^{a}\wedge \hat{\omega}_{a}{}^{\beta }, \\
& d\omega ^{n+\beta }=d\hat{\omega}^{n+\beta }=\sum_{a=1}^{2n}\hat{\omega}%
^{a}\wedge \hat{\omega}_{a}{}^{n+\beta }, \\
& d\hat{\omega}^{2n}=\sum_{a=1}^{2n}\hat{\omega}^{a}\wedge \hat{\omega}%
_{a}{}^{2n}.
\end{split}%
\end{equation}%
On the other hand, from the structure equations (see (\ref{A3})), we have 
\begin{equation}
\begin{split}
d\omega _{2n}{}^{\beta }& =\omega _{2n}{}^{\gamma }\wedge \omega _{\gamma
}{}^{\beta }+\omega _{2n}{}^{n+\gamma }\wedge \omega _{n+\gamma }{}^{\beta }(%
\text{summation convention)} \\
& =\sum_{\gamma =1}^{n-1}(-k\omega ^{\gamma }-\alpha \omega ^{n+\gamma
}+\Gamma ^{\gamma }\hat{\omega}^{2n})\wedge \omega _{\gamma }{}^{\beta
}+(-l\omega ^{n}+\Gamma ^{n}\hat{\omega}^{2n})\wedge \omega _{n}{}^{\beta }
\\
& \ \ \ +\sum_{\gamma =1}^{n-1}(\alpha \omega ^{\gamma }-k\omega ^{n+\gamma
}+\Gamma ^{n+\gamma }\hat{\omega}^{2n})\wedge \omega _{n+\gamma }{}^{\beta }.
\end{split}
\label{dofcon2}
\end{equation}%
We compare the coefficients of the terms $\hat{\omega}^{a}\wedge \hat{\omega}%
^{\beta }$ on both (\ref{dofcon1}) and (\ref{dofcon2}) for $n\geq 3$ and%
\newline
(i) $a=\gamma ,\ 1\leq \gamma \leq n-1$ : 
\begin{equation}
\begin{split}
& \text{The coefficient in}\ (\ref{dofcon1}) \\
& =-\hat{e}_{\gamma }k+k\hat{\omega}_{\beta }{}^{\beta }(\hat{e}_{\gamma })-k%
\hat{\omega}_{\gamma }{}^{\beta }(\hat{e}_{\beta })-\alpha \hat{\omega}%
_{\gamma }{}^{n+\beta }(\hat{e}_{\beta }) \\
& \ \ \ +\alpha \hat{\omega}_{\beta }{}^{n+\beta }(\hat{e}_{\gamma })+\Gamma
^{\beta }\hat{\omega}_{\gamma }{}^{2n}(\hat{e}_{\beta })-\Gamma ^{\beta }%
\hat{\omega}_{\beta }{}^{2n}(\hat{e}_{\gamma }) \\
& =-\hat{e}_{\gamma }k-k\omega _{\gamma }{}^{\beta }(e_{\beta })-\alpha
\omega _{\gamma }{}^{n+\beta }(e_{\beta })+\alpha \omega _{\beta
}{}^{n+\beta }(e_{\gamma }) \\
& \ \ \ +\frac{\alpha }{\sqrt{1+\alpha ^{2}}}\Gamma ^{\beta }\omega _{\gamma
}{}^{2n}(e_{\beta })-\frac{\alpha }{\sqrt{1+\alpha ^{2}}}\Gamma ^{\beta
}\omega _{\beta }{}^{2n}(e_{\gamma }) \\
& =-\hat{e}_{\gamma }k-k\omega _{\gamma }{}^{\beta }(e_{\beta })-\alpha
\omega _{\gamma }{}^{n+\beta }(e_{\beta })+\alpha \omega _{\beta
}{}^{n+\beta }(e_{\gamma }),
\end{split}%
\end{equation}%
and 
\begin{equation}
\begin{split}
& \text{The coefficient in}\ (\ref{dofcon2}) \\
& =-k\omega _{\gamma }{}^{\beta }(e_{\beta })+k\omega _{\beta }{}^{\beta
}(e_{\gamma })+\alpha \omega _{n+\gamma }{}^{\beta }(e_{\beta })-\alpha
\omega _{n+\beta }{}^{\beta }(e_{\gamma }) \\
& =-k\omega _{\gamma }{}^{\beta }(e_{\beta })+\alpha \omega _{n+\gamma
}{}^{\beta }(e_{\beta })-\alpha \omega _{n+\beta }{}^{\beta }(e_{\gamma }).
\end{split}%
\end{equation}%
Comparing the above two formulae, we get 
\begin{equation}
\hat{e}_{\gamma }k=0,\ \ 1\leq \gamma \leq n-1.  \label{parintcon1}
\end{equation}%
(ii) $a=n+\gamma ,\ 1\leq \gamma \leq n-1$ : 
\begin{equation}
\begin{split}
& \text{The coefficient in}\ (\ref{dofcon1}) \\
& =-\hat{e}_{n+\gamma }k+k\hat{\omega}_{\beta }{}^{\beta }(\hat{e}_{n+\gamma
})-k\hat{\omega}_{n+\gamma }{}^{\beta }(\hat{e}_{\beta })+(\hat{e}_{\beta
}\alpha )\delta _{\gamma \beta } \\
& \ \ \ -\alpha \hat{\omega}_{n+\gamma }{}^{n+\beta }(\hat{e}_{\beta
})+\alpha \hat{\omega}_{\beta }{}^{n+\beta }(\hat{e}_{n+\gamma })+\Gamma
^{\beta }\hat{\omega}_{n+\gamma }{}^{2n}(\hat{e}_{\beta })-\Gamma ^{\beta }%
\hat{\omega}_{\beta }{}^{2n}(\hat{e}_{n+\gamma }) \\
& =-\hat{e}_{n+\gamma }k-k\omega _{n+\gamma }{}^{\beta }(e_{\beta })+(\hat{e}%
_{\beta }\alpha )\delta _{\gamma \beta }-\alpha \omega _{n+\gamma
}{}^{n+\beta }(e_{\beta })+\alpha \omega _{\beta }{}^{n+\beta }(e_{n+\gamma
}) \\
& \ \ \ +\Gamma ^{\beta }\left( \frac{\alpha }{\sqrt{1+\alpha ^{2}}}\omega
_{n+\gamma }{}^{2n}(e_{\beta })-\frac{\delta _{\gamma \beta }}{\sqrt{%
1+\alpha ^{2}}}\right) -\Gamma ^{\beta }\left( \frac{\alpha }{\sqrt{1+\alpha
^{2}}}\omega _{\beta }{}^{2n}(e_{n+\gamma })+\frac{\delta _{\gamma \beta }}{%
\sqrt{1+\alpha ^{2}}}\right) ,
\end{split}%
\end{equation}%
and 
\begin{equation}
\text{The coefficient in}\ (\ref{dofcon2})=-k\omega _{n+\gamma }{}^{\beta
}(e_{\beta })+k\omega _{\beta }{}^{\beta }(e_{n+\gamma })-\alpha \omega
_{\gamma }{}^{\beta }(e_{\beta })-\alpha \omega _{n+\beta }{}^{\beta
}(e_{n+\gamma }).
\end{equation}%
Comparing the above two formulae, we get 
\begin{equation}
\hat{e}_{n+\gamma }k=0,\ \ 1\leq \gamma \leq n-1,\ \gamma \neq \beta ,
\label{parintcon2}
\end{equation}%
and 
\begin{equation}
\begin{split}
\hat{e}_{n+\beta }k& =\hat{e}_{\beta }\alpha +\Gamma ^{\beta }\left( \frac{%
\alpha }{\sqrt{1+\alpha ^{2}}}(h_{(n+\beta )\beta }-h_{\beta (n+\beta )})-%
\frac{2}{\sqrt{1+\alpha ^{2}}}\right)  \\
& =3\hat{e}_{\beta }\alpha .
\end{split}
\label{parintcon3}
\end{equation}%
Since $\beta $ is arbitrary, from (\ref{parintcon2}) and (\ref{parintcon3}),
we conclude that 
\begin{equation}
\hat{e}_{n+\gamma }k=\hat{e}_{\gamma }\alpha =0,\ \ 1\leq \gamma \leq n-1.%
\text{ }n\geq 3  \label{parintcon4}
\end{equation}%
Moreover, we have 
\begin{equation}
\Gamma ^{\gamma }=\omega _{2n}{}^{\gamma }(\hat{e}_{2n})=-\left( \frac{\hat{e%
}_{\gamma }\alpha }{\sqrt{1+\alpha ^{2}}}\right) =0.  \label{parintcon5}
\end{equation}%
\noindent The second equality in (\ref{parintcon5}) is due to (\ref{2-0}).
Notice that, for $n=2,$ we only get formula (\ref{parintcon3}). Similarly,
fixing $\beta ,1\leq \beta \leq n-1$, if we take the exterior differential
of the second formula of (\ref{expofcon}) and compare the coefficients of
the terms $\hat{\omega}^{a}\wedge \hat{\omega}^{n+\beta }$ for $1\leq a\leq
2n-1,\ a\neq n$, we get 
\begin{equation}
\hat{e}_{\gamma }k=0,\ \ 1\leq \gamma \leq n-1,\text{ }n\geq 3,
\label{parintcon6}
\end{equation}%
\begin{equation}
\hat{e}_{\gamma }k=-3\hat{e}_{n+\gamma }\alpha .  \label{0.17}
\end{equation}%
and hence 
\begin{equation}
\Gamma ^{n+\gamma }=\omega _{2n}{}^{n+\gamma }(\hat{e}_{2n})=-\left( \frac{%
\hat{e}_{n+\gamma }\alpha }{\sqrt{1+\alpha ^{2}}}\right) =0.
\label{parintcon7}
\end{equation}%
\noindent The second equality in (\ref{parintcon7}) is due to (\ref{2-0}).
Again, for $n=2,$ we just get formula (\ref{0.17}). In order to show the
formula $0=\hat{e}_{\gamma }k=\hat{e}_{n+\gamma }k$, or equivalently to show
that $\hat{e}_{\gamma }\alpha =\hat{e}_{n+\gamma }\alpha =0$ for $n\geq 2$,
we need to take the exterior differential of the third equation of (\ref%
{expofcon}).

From the structure equation, we have 
\begin{equation}
\begin{split}
d\omega _{n}{}^{2n}& =\sum_{\gamma =1}^{n}\omega _{n}{}^{\gamma }\wedge
\omega _{\gamma }{}^{2n}+\omega _{n}{}^{n+\gamma }\wedge \omega _{n+\gamma
}{}^{2n} \\
& =\sum_{\gamma =1}^{n}\omega _{2n}{}^{n+\gamma }\wedge \omega _{\gamma
}{}^{2n}-\omega _{2n}{}^{\gamma }\wedge \omega _{n+\gamma }{}^{2n} \\
& =2\sum_{\gamma =1}^{n-1}\omega _{\gamma }{}^{2n}\wedge \omega _{n+\gamma
}{}^{2n} \\
& =2\sum_{\gamma =1}^{n-1}(k\omega ^{\gamma }+\alpha \omega ^{n+\gamma
}-\Gamma ^{\gamma }\hat{\omega}^{2n})\wedge (-\alpha \omega ^{\gamma
}+k\omega ^{n+\gamma }-\Gamma ^{n+\gamma }\hat{\omega}^{2n}) \\
& =2\left( \sum_{\gamma =1}^{n-1}(k^{2}+\alpha ^{2})\omega ^{\gamma }\wedge
\omega ^{n+\gamma }-(\alpha \Gamma ^{\gamma }+k\Gamma ^{n+\gamma })\omega
^{\gamma }\wedge \hat{\omega}^{2n}+(k\Gamma ^{\gamma }-\alpha \Gamma
^{n+\gamma })\omega ^{n+\gamma }\wedge \hat{\omega}^{2n}\right) .
\end{split}
\label{dofcon3}
\end{equation}%
On the other hand, taking the exterior differential of the third equation of
(\ref{expofcon}), we have 
\begin{equation}
d\omega _{n}{}^{2n}=\sum_{a=1}^{2n}(\hat{e}_{a}l)\hat{\omega}^{a}\wedge
\omega ^{n}+ld\omega ^{n}-\sum_{a=1}^{2n}(\hat{e}_{a}\Gamma ^{n})\hat{\omega}%
^{a}\wedge \hat{\omega}^{2n}-\Gamma ^{n}d\hat{\omega}^{2n},  \label{dofcon4}
\end{equation}%
where 
\begin{equation}
\begin{split}
d\omega ^{n}& =d\hat{\omega}^{n}=\sum_{a=1}^{2n}\hat{\omega}^{a}\wedge \hat{%
\omega}_{a}{}^{n} \\
& =\left( \sum_{\gamma =1}^{n-1}\hat{\omega}^{\gamma }\wedge \hat{\omega}%
_{\gamma }{}^{n}+\hat{\omega}^{n+\gamma }\wedge \hat{\omega}_{n+\gamma
}{}^{n}\right) +\hat{\omega}^{2n}\wedge \hat{\omega}_{2n}{}^{n} \\
& =\sum_{\gamma =1}^{n-1}\hat{\omega}^{\gamma }\wedge (-\alpha \omega
^{\gamma }+k\omega ^{n+\gamma }-\Gamma ^{n+\gamma }\hat{\omega}^{2n})+\hat{%
\omega}^{n+\gamma }(-k\omega ^{\gamma }-\alpha \omega ^{n+\gamma }+\Gamma
^{\gamma }\hat{\omega}^{2n}) \\
& \ \ \ +\frac{\alpha }{\sqrt{1+\alpha ^{2}}}\hat{\omega}^{2n}\wedge
(-l\omega ^{n}+\Gamma ^{n}\hat{\omega}^{2n}) \\
& =\left( \sum_{\gamma =1}^{n-1}2k\hat{\omega}^{\gamma }\wedge \hat{\omega}%
^{n+\gamma }-\Gamma ^{n+\gamma }\hat{\omega}^{\gamma }\wedge \hat{\omega}%
^{2n}+\Gamma ^{\gamma }\hat{\omega}^{n+\gamma }\wedge \hat{\omega}%
^{2n}\right) +\frac{l\alpha }{\sqrt{1+\alpha ^{2}}}\hat{\omega}^{n}\wedge 
\hat{\omega}^{2n},
\end{split}
\label{parintcon17}
\end{equation}%
and 
\begin{equation}
\begin{split}
d\hat{\omega}^{2n}& =\sum_{a=1}^{2n}\hat{\omega}^{a}\wedge \hat{\omega}%
_{a}{}^{2n} \\
& =\left( \sum_{\gamma =1}^{n-1}\hat{\omega}^{\gamma }\wedge \hat{\omega}%
_{\gamma }{}^{2n}+\hat{\omega}^{n+\gamma }\wedge \hat{\omega}_{n+\gamma
}{}^{2n}\right) +\hat{\omega}^{n}\wedge \hat{\omega}_{n}{}^{2n} \\
& =\sum_{\gamma =1}^{n-1}\hat{\omega}^{\gamma }\wedge \left( \frac{\alpha }{%
\sqrt{1+\alpha ^{2}}}(k\omega ^{\gamma }+\alpha \omega ^{n+\gamma }-\Gamma
^{\gamma }\hat{\omega}^{2n})+\frac{1}{\sqrt{1+\alpha ^{2}}}\hat{\omega}%
^{n+\gamma }\right) \\
& \ \ \ +\sum_{\gamma =1}^{n-1}\hat{\omega}^{n+\gamma }\wedge \left( \frac{%
\alpha }{\sqrt{1+\alpha ^{2}}}(-\alpha \omega ^{\gamma }+k\omega ^{n+\gamma
}-\Gamma ^{n+\gamma }\hat{\omega}^{2n})-\frac{1}{\sqrt{1+\alpha ^{2}}}\hat{%
\omega}^{\gamma }\right) \\
& \ \ \ +\hat{\omega}^{n}\wedge \left( \frac{\alpha }{\sqrt{1+\alpha ^{2}}}%
(l\omega ^{n}-\Gamma ^{n}\hat{\omega}^{2n})+\frac{2\alpha }{1+\alpha ^{2}}%
\hat{\omega}^{2n}\right) \\
& =\left( \sum_{\gamma =1}^{n-1}2\sqrt{1+\alpha ^{2}}\hat{\omega}^{\gamma
}\wedge \hat{\omega}^{n+\gamma }-\frac{\alpha }{\sqrt{1+\alpha ^{2}}}\Gamma
^{\gamma }\hat{\omega}^{\gamma }\wedge \hat{\omega}^{2n}-\frac{\alpha }{%
\sqrt{1+\alpha ^{2}}}\Gamma ^{n+\gamma }\hat{\omega}^{n+\gamma }\wedge \hat{%
\omega}^{2n}\right) \\
& \ \ \ +\left( \frac{2\alpha }{1+\alpha ^{2}}-\frac{\alpha }{\sqrt{1+\alpha
^{2}}}\Gamma ^{n}\right) \hat{\omega}^{n}\wedge \hat{\omega}^{2n}.
\end{split}
\label{parintcon18}
\end{equation}%
Substituting (\ref{parintcon17}) and (\ref{parintcon18}) into (\ref{dofcon4}%
), we get 
\begin{equation}
\begin{split}
d\omega _{n}{}^{2n}& =\sum_{\gamma =1}^{n-1}(\hat{e}_{\gamma }l)\hat{\omega}%
^{\gamma }\wedge \hat{\omega}^{n}-\sum_{\gamma =1}^{n-1}(\hat{e}_{n+\gamma
}l)\hat{\omega}^{n}\wedge \hat{\omega}^{n+\gamma }-\sum_{a=1}^{2n}(\hat{e}%
_{a}\Gamma ^{n})\hat{\omega}^{a}\wedge \hat{\omega}^{2n}-\hat{e}_{2n}l\hat{%
\omega}^{n}\wedge \hat{\omega}^{2n} \\
& \ \ \ +(\sum_{\gamma =1}^{n-1}2kl\hat{\omega}^{\gamma }\wedge \hat{\omega}%
^{n+\gamma })+\frac{l^{2}\alpha }{\sqrt{1+\alpha ^{2}}}\hat{\omega}%
^{n}\wedge \hat{\omega}^{2n} \\
& \ \ \ -\left( \sum_{\gamma =1}^{n-1}2\sqrt{1+\alpha ^{2}}\Gamma ^{n}\hat{%
\omega}^{\gamma }\wedge \hat{\omega}^{n+\gamma }\right) -\left( \frac{%
2\alpha }{1+\alpha ^{2}}-\frac{\alpha }{\sqrt{1+\alpha ^{2}}}\Gamma
^{n}\right) \Gamma ^{n}\hat{\omega}^{n}\wedge \hat{\omega}^{2n} \\
& \ \ \ +\left( \frac{\alpha }{\sqrt{1+\alpha ^{2}}}\Gamma ^{\gamma }\Gamma
^{n}-l\Gamma ^{n+\gamma }\right) \hat{\omega}^{\gamma }\wedge \hat{\omega}%
^{2n}+\left( \frac{\alpha }{\sqrt{1+\alpha ^{2}}}\Gamma ^{n+\gamma }\Gamma
^{n}+l\Gamma ^{\gamma }\right) \hat{\omega}^{n+\gamma }\wedge \hat{\omega}%
^{2n}.
\end{split}
\label{dofcon5}
\end{equation}%
Comparing, respectively, the coefficients of both the term $\hat{\omega}%
^{\gamma }\wedge \hat{\omega}^{n+\gamma }$ and $\hat{\omega}^{a}\wedge \hat{%
\omega}^{n},a\neq n,2n$ of (\ref{dofcon3}) and (\ref{dofcon5}), we get 
\begin{equation}
\begin{split}
\hat{e}_{n}\alpha & =k^{2}-\alpha ^{2}-kl \\
\hat{e}_{\gamma }l& =\hat{e}_{n+\gamma }l=0,\ \ \text{for}\ 1\leq \gamma
\leq n-1.
\end{split}
\label{n2f4}
\end{equation}%
For $n=2$, if we compare the coefficients of both the term $\hat{\omega}%
^{1}\wedge \hat{\omega}^{4}$ and $\hat{\omega}^{3}\wedge \hat{\omega}^{4}$
of (\ref{dofcon3}) and (\ref{dofcon5}), we get 
\begin{equation}
\begin{split}
-2\alpha \Gamma ^{1}-2k\Gamma ^{3}& =-\hat{e}_{1}\Gamma ^{2}+\frac{\alpha }{%
\sqrt{1+\alpha ^{2}}}\Gamma ^{1}\Gamma ^{2}-l\Gamma ^{3} \\
2k\Gamma ^{1}-2\alpha \Gamma ^{3}& =-\hat{e}_{3}\Gamma ^{2}+\frac{\alpha }{%
\sqrt{1+\alpha ^{2}}}\Gamma ^{3}\Gamma ^{2}+l\Gamma ^{1},
\end{split}
\label{n2f3}
\end{equation}%
where, from (\ref{2-0}) and (\ref{n2f4}), we have 
\begin{equation*}
\Gamma ^{2}=\frac{-2\alpha ^{2}-\hat{e}_{2}\alpha }{\sqrt{1+\alpha ^{2}}}=%
\frac{kl-\alpha ^{2}-k^{2}}{\sqrt{1+\alpha ^{2}}},
\end{equation*}%
hence 
\begin{equation}
\hat{e}_{a}\Gamma ^{2}=\frac{(1+\alpha ^{2})(l-2k)\hat{e}_{a}k-(2\alpha
+\alpha kl-\alpha ^{3}-\alpha k^{2})\hat{e}_{a}\alpha }{(1+\alpha ^{2})^{3/2}%
},\ \ \text{for}\ a=1,3.  \label{n2f5}
\end{equation}%
Substituting (\ref{n2f5}) into (\ref{n2f3}), using formulae (\ref{parintcon3}%
), (\ref{0.17}), and noting that $\Gamma ^{a}=\frac{-\hat{e}_{a}\alpha }{%
\sqrt{1+\alpha ^{2}}},\ a=1,3$ by (\ref{2-0}), we get 
\begin{equation}
\begin{array}{rl}
\alpha ^{3}\hat{e}_{1}\alpha -2(1+\alpha ^{2})(l-2k)\hat{e}_{3}\alpha & =0
\\ 
2(1+\alpha ^{2})(l-2k)\hat{e}_{1}\alpha +\alpha ^{3}\hat{e}_{3}\alpha & =0.%
\end{array}
\label{n2f6}
\end{equation}%
It is easy to see that the determinant of the coefficients matrix of
equations (\ref{n2f6}) is $\alpha ^{6}+4(1+\alpha ^{2})^{2}(l-2k)^{2}$,
hence it vanishes if and only if $\alpha =0$ and $l=2k$. Together with (\ref%
{n2f4}), we get $\alpha =l=k=0$, which are all zero. If the determinant of
the coefficients matrix is not zero, then we immediately have $\hat{e}%
_{1}\alpha =\hat{e}_{3}\alpha =0$, thus also $\hat{e}_{1}k=\hat{e}%
_{3}k=\Gamma ^{1}=\Gamma ^{3}=0$, for $n=2$.

Now we continue to compare the coefficients of the terms $\hat{\omega}%
^{a}\wedge \hat{\omega}^{\beta }$ on both (\ref{dofcon1}) and (\ref{dofcon2}%
) for $a=n$. We have 
\begin{equation}
\begin{split}
& \text{The coefficient in}\ (\ref{dofcon1}) \\
& =-\hat{e}_{n}k-k\hat{\omega}_{n}{}^{\beta }(\hat{e}_{\beta })+k\hat{\omega}%
_{\beta }{}^{\beta }(\hat{e}_{n})-\alpha \hat{\omega}_{n}{}^{n+\beta }(\hat{e%
}_{\beta }) \\
& \ \ \ +\alpha \hat{\omega}_{\beta }{}^{n+\beta }(\hat{e}_{n})+\Gamma
^{\beta }\hat{\omega}_{n}{}^{2n}(\hat{e}_{\beta })-\Gamma ^{\beta }\hat{%
\omega}_{\beta }{}^{2n}(\hat{e}_{n}) \\
& =-\hat{e}_{n}k-k\alpha -\alpha k+\alpha \omega _{\beta }{}^{n+\beta
}(e_{n}),
\end{split}
\label{parintcon8}
\end{equation}%
where, for the last equality, we have used 
\begin{equation}
\begin{split}
\hat{\omega}_{n}{}^{\beta }(\hat{e}_{\beta })& =\omega _{n}{}^{\beta
}(e_{\beta })=-\omega _{n+\beta }{}^{2n}(e_{\beta })=-h_{(n+\beta )\beta
}=\alpha \\
\hat{\omega}_{n}{}^{n+\beta }(\hat{e}_{\beta })& =\omega _{n}{}^{n+\beta
}(e_{\beta })=\omega _{\beta }{}^{2n}(e_{\beta })=h_{\beta \beta }=k,
\end{split}%
\end{equation}%
and 
\begin{equation}
\begin{split}
\text{The coefficient in}\ (\ref{dofcon2})& =k\omega _{\beta }{}^{\beta }(%
\hat{e}_{n})-l\omega _{n}{}^{\beta }(\hat{e}_{\beta })-\alpha \omega
_{n+\beta }{}^{\beta }(\hat{e}_{n}) \\
& =-l\alpha -\alpha \omega _{n+\beta }{}^{\beta }(\hat{e}_{n}).
\end{split}
\label{parintcon9}
\end{equation}%
From (\ref{parintcon8}) and (\ref{parintcon9}), we get 
\begin{equation}
\hat{e}_{n}k=(l-2k)\alpha .  \label{parintcon10}
\end{equation}%
Then, we compare the coefficients of the terms $\hat{\omega}^{2n}\wedge \hat{%
\omega}^{\beta }$ on both (\ref{dofcon1}) and (\ref{dofcon2}). We have 
\begin{equation}
\begin{split}
& \text{The coefficient in}\ (\ref{dofcon1}) \\
& =-\hat{e}_{2n}k-k\hat{\omega}_{2n}{}^{\beta }(\hat{e}_{\beta })+k\hat{%
\omega}_{\beta }{}^{\beta }(\hat{e}_{2n})-\alpha \hat{\omega}%
_{2n}{}^{n+\beta }(\hat{e}_{\beta }) \\
& \ \ \ +\alpha \hat{\omega}_{\beta }{}^{n+\beta }(\hat{e}_{2n})-(\hat{e}%
_{\beta }\Gamma ^{\beta })+\Gamma ^{\beta }\hat{\omega}_{2n}{}^{2n}(\hat{e}%
_{\beta })-\Gamma ^{\beta }\hat{\omega}_{\beta }{}^{2n}(\hat{e}_{2n}) \\
& =-\hat{e}_{2n}k+\frac{k^{2}\alpha }{\sqrt{1+\alpha ^{2}}}-\alpha \sqrt{%
1+\alpha ^{2}}+\alpha \hat{\omega}_{\beta }{}^{n+\beta }(\hat{e}_{2n}),
\end{split}
\label{parintcon11}
\end{equation}%
where, for the last equality, we have used 
\begin{equation}
\begin{split}
\hat{\omega}_{2n}{}^{\beta }(\hat{e}_{\beta })& =-\frac{\alpha }{\sqrt{%
1+\alpha ^{2}}}\omega _{\beta }{}^{2n}(e_{\beta })-\frac{1}{\sqrt{1+\alpha
^{2}}}\hat{\omega}^{n+\beta }(e_{\beta }) \\
& =-\frac{\alpha }{\sqrt{1+\alpha ^{2}}}h_{\beta \beta }=-\frac{k\alpha }{%
\sqrt{1+\alpha ^{2}}},
\end{split}%
\end{equation}%
and 
\begin{equation}
\begin{split}
\hat{\omega}_{2n}{}^{n+\beta }(\hat{e}_{\beta })& =-\frac{\alpha }{\sqrt{%
1+\alpha ^{2}}}\omega _{n+\beta }{}^{2n}(e_{\beta })+\frac{1}{\sqrt{1+\alpha
^{2}}}\hat{\omega}^{\beta }(e_{\beta }) \\
& =-\frac{\alpha }{\sqrt{1+\alpha ^{2}}}h_{(n+\beta )\beta }+\frac{1}{\sqrt{%
1+\alpha ^{2}}}=\sqrt{1+\alpha ^{2}},
\end{split}%
\end{equation}%
and 
\begin{equation}
\begin{split}
\text{The coefficient in}\ (\ref{dofcon2})& =k\omega _{\beta }{}^{\beta }(%
\hat{e}_{2n})-\alpha \omega _{n+\beta }{}^{\beta }(\hat{e}%
_{2n})+\sum_{j=1}^{2n-1}\Gamma ^{j}\omega _{j}{}^{\beta }(\hat{e}_{\beta })
\\
& =-\alpha \hat{\omega}_{n+\beta }{}^{\beta }(\hat{e}_{2n})-\frac{\alpha }{%
\sqrt{1+\alpha ^{2}}}+\Gamma ^{n}\omega _{n}{}^{\beta }(\hat{e}_{\beta }) \\
& =-\alpha \hat{\omega}_{n+\beta }{}^{\beta }(\hat{e}_{2n})-\frac{\alpha }{%
\sqrt{1+\alpha ^{2}}}-\left( \frac{\hat{e}_{n}\alpha +2\alpha ^{2}}{\sqrt{%
1+\alpha ^{2}}}\right) \omega _{2n}{}^{n+\beta }(e_{\beta }) \\
& =-\alpha \hat{\omega}_{n+\beta }{}^{\beta }(\hat{e}_{2n})-\frac{\alpha }{%
\sqrt{1+\alpha ^{2}}}+\left( \frac{\hat{e}_{n}\alpha +2\alpha ^{2}}{\sqrt{%
1+\alpha ^{2}}}\right) h_{(n+\beta )\beta } \\
& =-\alpha \hat{\omega}_{n+\beta }{}^{\beta }(\hat{e}_{2n})-\frac{\alpha }{%
\sqrt{1+\alpha ^{2}}}-\left( \frac{\alpha (\hat{e}_{n}\alpha +2\alpha ^{2})}{%
\sqrt{1+\alpha ^{2}}}\right) .
\end{split}
\label{parintcon12}
\end{equation}%
From (\ref{parintcon11}) and (\ref{parintcon12}), we get 
\begin{equation}
\hat{e}_{2n}k=\frac{\alpha }{\sqrt{1+\alpha ^{2}}}(k^{2}+\alpha ^{2}+\hat{e}%
_{n}\alpha ).  \label{parintcon13}
\end{equation}%
Finally, we compare the coefficients of the terms $\hat{\omega}^{2n}\wedge 
\hat{\omega}^{n+\beta }$ on both (\ref{dofcon1}) and (\ref{dofcon2}). We
have 
\begin{equation}
\begin{split}
\text{The coefficient in}\ (\ref{dofcon1})& =-k\hat{\omega}_{2n}{}^{\beta }(%
\hat{e}_{n+\beta })+k\hat{\omega}_{n+\beta }{}^{\beta }(\hat{e}_{2n}) \\
& \ \ \ -(\hat{e}_{2n}\alpha )-\alpha \hat{\omega}_{2n}{}^{n+\beta }(\hat{e}%
_{n+\beta })+\alpha \hat{\omega}_{n+\beta }{}^{n+\beta }(\hat{e}_{2n}) \\
& =k\sqrt{1+\alpha ^{2}}+k\hat{\omega}_{n+\beta }{}^{\beta }(\hat{e}_{2n})-(%
\hat{e}_{2n}\alpha )+\frac{\alpha ^{2}k}{\sqrt{1+\alpha ^{2}}},
\end{split}
\label{parintcon14}
\end{equation}%
where, for the last equality, we have used 
\begin{equation}
\begin{split}
\hat{\omega}_{2n}{}^{n+\beta }(\hat{e}_{n+\beta })& =-\frac{\alpha }{\sqrt{%
1+\alpha ^{2}}}\omega _{n+\beta }{}^{2n}(e_{n+\beta })+\frac{1}{\sqrt{%
1+\alpha ^{2}}}\hat{\omega}^{\beta }(e_{n+\beta }) \\
& =-\frac{\alpha }{\sqrt{1+\alpha ^{2}}}h_{(n+\beta )(n+\beta )}=-\frac{%
k\alpha }{\sqrt{1+\alpha ^{2}}},
\end{split}%
\end{equation}%
and 
\begin{equation}
\begin{split}
\hat{\omega}_{2n}{}^{\beta }(\hat{e}_{n+\beta })& =-\frac{\alpha }{\sqrt{%
1+\alpha ^{2}}}\omega _{\beta }{}^{2n}(e_{n+\beta })-\frac{1}{\sqrt{1+\alpha
^{2}}}\hat{\omega}^{n+\beta }(\hat{e}_{n+\beta }) \\
& =-\frac{\alpha }{\sqrt{1+\alpha ^{2}}}h_{\beta (n+\beta )}-\frac{1}{\sqrt{%
1+\alpha ^{2}}}=-\sqrt{1+\alpha ^{2}},
\end{split}%
\end{equation}%
and 
\begin{equation}
\begin{split}
\text{The coefficient in}\ (\ref{dofcon2})& =\alpha \omega _{\beta
}{}^{\beta }(\hat{e}_{2n})+k\omega _{n+\beta }{}^{\beta }(\hat{e}%
_{2n}+\Gamma ^{n}\omega _{n}{}^{\beta }(\hat{e}_{n+\beta }) \\
& =-k\left( \hat{\omega}_{\beta }{}^{n+\beta }(\hat{e}_{2n})-\frac{1}{\sqrt{%
1+\alpha ^{2}}}\right) +\Gamma ^{n}\omega _{2n}{}^{n+\beta }(\hat{e}%
_{n+\beta }) \\
& =-k\hat{\omega}_{\beta }{}^{n+\beta }(\hat{e}_{2n})+\frac{k}{\sqrt{%
1+\alpha ^{2}}}-k\Gamma ^{n}.
\end{split}
\label{parintcon15}
\end{equation}%
From (\ref{parintcon14}) and (\ref{parintcon15}), we get 
\begin{equation}
\begin{split}
\hat{e}_{2n}\alpha & =k\sqrt{1+\alpha ^{2}}+\frac{k\alpha ^{2}}{\sqrt{%
1+\alpha ^{2}}}-\frac{k}{\sqrt{1+\alpha ^{2}}}+k\Gamma ^{n} \\
& =\frac{2k\alpha ^{2}}{\sqrt{1+\alpha ^{2}}}-k\left( \frac{\hat{e}%
_{n}\alpha +2\alpha ^{2}}{\sqrt{1+\alpha ^{2}}}\right) =-k\frac{\hat{e}%
_{n}\alpha }{\sqrt{1+\alpha ^{2}}}.
\end{split}
\label{parintcon16}
\end{equation}

Notice that we have shown $\Gamma ^{a}=0$, for $1\leq a\leq 2n-1,\ a\neq n$,
so if we again compare the coefficients of (\ref{dofcon3}) and (\ref{dofcon5}%
), we have 
\begin{equation}
\hat{e}_{a}\Gamma ^{n}=0,\ \ \text{for}\ 1\leq a\leq 2n,\ a\neq n,2n,
\label{parintcon19}
\end{equation}%
and 
\begin{equation}
0=\frac{l^{2}\alpha }{\sqrt{1+\alpha ^{2}}}-\left( \frac{2\alpha }{1+\alpha
^{2}}-\frac{\alpha }{\sqrt{1+\alpha ^{2}}}\Gamma ^{n}\right) \Gamma ^{n}-%
\hat{e}_{2n}l-\hat{e}_{n}\Gamma ^{n}.  \label{parintcon20}
\end{equation}%
Since $\hat{e}_{a}\alpha =0$, we have 
\begin{equation}
\hat{e}_{a}\Gamma ^{n}=-\frac{\hat{e}_{a}\hat{e}_{n}\alpha }{\sqrt{1+\alpha
^{2}}}
\end{equation}%
\noindent by (\ref{2-0}). Observe that (\ref{parintcon19}) is equivalentt to 
\begin{equation}
\hat{e}_{a}\hat{e}_{n}\alpha =0,\ \ \text{for}\ 1\leq a\leq 2n,\ a\neq n,2n.
\end{equation}%
After a direct computation, we see that (\ref{parintcon20}) is just a
Codazzi-like equation, which is the last equation of (\ref{parintcon}).
Therefore we have completed the proof of Proposition 4.2.

\bigskip

\section{An ODE system and proof of Lemma B}

From Proposition 4.2, $k$ and $\alpha $ satisfies the following equations 
\begin{eqnarray}
e_{n}k &=&(l-2k)\alpha  \label{5.1} \\
e_{n}\alpha &=&k^{2}-\alpha ^{2}-kl  \notag
\end{eqnarray}%
\noindent on an umbilic hypersurface $\Sigma $ of $H_{n}$. Observe that $p$%
-mean curvature $H$ of $\Sigma $ and $k,$ $l$ have the following relation: 
\begin{equation}
H=(2n-2)k+l.  \label{5.2}
\end{equation}%
\noindent Let $\beta $ :$=$ $l-2k$ and write $e_{n}k,$ $e_{n}\alpha $, etc.
as $k^{\prime }$, $\alpha ^{\prime }$, etc.. We can then express (\ref{5.1})
in terms of $\beta $, $\alpha $ as: 
\begin{eqnarray}
\beta ^{\prime } &=&-2n\beta \alpha  \label{5.3} \\
\alpha ^{\prime } &=&-\alpha ^{2}+\frac{1}{4n^{2}}(\beta -c)((2n-1)\beta +c)
\notag
\end{eqnarray}%
\noindent on $\Sigma $ having $H$ $=$ $c$, a positive constant, by (\ref{5.2}%
). Let $\Lambda $ denote the set in the $\alpha \beta $-plane, which
consists of 
\begin{equation}
\beta =0\text{ }(\alpha \text{-}axis)  \label{5.4}
\end{equation}%
\noindent (which is a solution to (\ref{5.3}) with $\alpha ^{\prime }$ $=$ $%
-\alpha ^{2}$ $-$ $\frac{c^{2}}{4n^{2}}$ $<$ $0)$ and two points:%
\begin{equation}
\alpha =0,\text{ }\beta =c\text{ or }-\frac{c}{2n-1}  \label{5.5}
\end{equation}

\noindent (which are stationary points of (\ref{5.3})). Write $%
R^{2}\backslash \Lambda $ $=$ $R^{2,+}\backslash \{(0,c)\}$ $\cup $ $%
R^{2,-}\backslash \{(0,-\frac{c}{2n-1})\}$ where $R^{2,+}$ ($R^{2,-},$
resp.) $:=$ $\{\beta $ $>$ $0\}$ ($\{\beta $ $<$ $0\},$ resp.)$.$

$\bigskip $

\textbf{Lemma 6.1. }\textit{For any initial point }$p_{0}$\textit{\ }$=$%
\textit{\ }$(\alpha _{0},\beta _{0})$\textit{\ }$\in $\textit{\ }$%
R^{2,+}\backslash \{(0,c)\}$ \textit{(}$R^{2,-}\backslash \{(0,-\frac{c}{2n-1%
})\},$\textit{\ resp.)}$,$\textit{\ there passes a unique periodic orbit }$%
\gamma $\textit{\ }$\subset $ $R^{2,+}\backslash \{(0,c)\}$ \textit{(}$%
R^{2,-}\backslash \{(0,-\frac{c}{2n-1})\}\mathit{,}$\textit{\ resp.),} 
\textit{described by }$(\alpha (s),\beta (s))$\textit{, }$0$\textit{\ }$\leq 
$\textit{\ }$s$\textit{\ }$\leq $\textit{\ }$s_{0},$\textit{\ which is a
solution to the }$ODE$\textit{\ system (\ref{5.3}), with }$\alpha (s_{0})$%
\textit{\ }$=$\textit{\ }$\alpha (0)$\textit{\ }$=$\textit{\ }$\alpha _{0}$%
\textit{\ and }$\beta (s_{0})$\textit{\ }$=$\textit{\ }$\beta (0)$\textit{\ }%
$=$\textit{\ }$\beta _{0}.$\textit{\ Moreover, }$\gamma $\textit{\ is
symmetric with respect to the }$\beta $\textit{-axis, i.e., }$(\alpha ,\beta
)$\textit{\ }$\in $\textit{\ }$\gamma $\textit{\ implies }$(-\alpha ,\beta )$%
\textit{\ }$\in $\textit{\ }$\gamma .$

\textit{\bigskip }

\proof
(I) Suppose $p_{0}$\textit{\ }$=$\textit{\ }$(\alpha _{0},\beta _{0})$%
\textit{\ }$\in $\textit{\ }$R^{2,+}\backslash \{(0,c)\}.$ Let $\Upsilon $
denote the hyperbolic curve in the $\alpha \beta $-plane defined by 
\begin{equation*}
(\alpha ^{\prime }=)-\alpha ^{2}+\frac{1}{4n^{2}}(\beta -c)((2n-1)\beta
+c)=0.
\end{equation*}%
\noindent Note that $\Upsilon $ passes through two (stationary) points $%
(0,c) $ and $(0,-\frac{c}{2n-1})$ (cf. (5.5)). Observe that $\Upsilon $ is
invariant under the reflection $(\alpha ,\beta )$ $\rightarrow $ $(-\alpha
,\beta )$ with respect to the $\beta $-axis, and equation (\ref{5.3}) has
the symmetry property that if ($\alpha ^{\prime }$, $\beta ^{\prime }$) at $%
(\alpha _{1},$ $\beta _{1})$ satisfies (\ref{5.3}), then $(\alpha ^{\prime
},-\beta ^{\prime })$ at $(-\alpha _{1},$ $\beta _{1})$ also satisfies (\ref%
{5.3}). So without loss of generality, we may assume $p_{0}$\textit{\ }$=$%
\textit{\ }$(\alpha _{0},\beta _{0})$ lies in the right half plane. Note
that $\Upsilon $ divides the first quadrant into two regions: 
\begin{eqnarray}
R_{+} &:&=\{(\alpha ,\beta ):\alpha >0,\beta >0,\alpha ^{\prime }>0\}
\label{5.5.1} \\
R_{-} &:&=\{(\alpha ,\beta ):\alpha >0,\beta >0,\alpha ^{\prime }<0\}. 
\notag
\end{eqnarray}

Let $V$ $:=$ $V(\alpha ,\beta )$ denote the following vector field at $%
(\alpha ,\beta ):$%
\begin{equation}
(-\alpha ^{2}+\frac{1}{4n^{2}}(\beta -c)((2n-1)\beta +c),-2n\beta \alpha ).
\label{5.6}
\end{equation}

\bigskip

\textbf{Case 1.} $p_{-}$ $=$\textit{\ }$(\alpha _{-},\beta _{-})$ $\in $ $%
\beta $-axis (hence $\alpha _{-}$ $=$ $0$) with $\beta _{-}$ $>$ $c$ $(>$ $%
0).$ Then there is small $\varepsilon $ $>$ $0$ such that the solution $%
p(s):=(\alpha (s),\beta (s))$ to (\ref{5.3}) with $p(s_{-})=p_{-}$ enters $%
R_{+}$ (the second quadrant, resp.) for $s_{-}$ $<$ $s$ $<$ $s_{-}$ $+$ $%
\varepsilon $ ($s_{-}$ $-$ $\varepsilon $ $<$ $s$ $<$ $s_{-},$ resp.) since 
\begin{equation*}
\alpha ^{\prime }=\frac{1}{4n^{2}}(\beta _{0}-c)((2n-1)\beta _{0}+c)>0
\end{equation*}%
\noindent ($\beta ^{\prime }$ $=$ $0)$ at $p_{-}$.

\textbf{Case 2.} $p_{0}$\textit{\ }$=$\textit{\ }$(\alpha _{0},\beta _{0})$ $%
\in $ $R_{+}$ (see (\ref{5.5.1})). Let $p(s):=(\alpha (s),\beta (s))$ denote
the solution to (\ref{5.3}) with $p(s_{0})$ $=$ $p_{0}$. Since $\alpha
^{\prime }$ $>$ $0$ and $\beta ^{\prime }$ $<$ $0$ in $R_{+},$ $\alpha $ is
decreasing while $\beta $ is increasing as time changes towards negative
infinity. Observe that%
\begin{eqnarray*}
\alpha ^{\prime } &=&-\alpha ^{2}+\frac{1}{4n^{2}}(\beta -c)((2n-1)\beta +c)
\\
&\geq &-\alpha _{0}^{2}+\frac{1}{4n^{2}}(\beta _{0}-c)((2n-1)\beta _{0}+c)>0
\end{eqnarray*}%
\noindent for $s$ $\leq $ $s_{0}.$ Therefore at a finite time $\tilde{s}_{0}$
$<$ $s_{0}$, $\alpha (\tilde{s}_{0})$ $=$ $0,$ i.e., $p(\tilde{s}_{0})$ $\in 
$ $\beta $-axis. On the other hand, as time changes towards the positive
infinity, $\alpha $ is increasing while $\beta $ is decreasing. Moreover, we
observe that%
\begin{equation*}
\beta ^{\prime }=-2n\beta \alpha \leq -2nc\alpha _{0}<0
\end{equation*}%
\noindent since $\beta $ $\geq $ $c$ in $R_{+}.$ Therefore $p(s)$ must hit $%
\Upsilon \backslash \{(0,c)\}$ $\cap $ (first quadrant) at a finite time $%
\breve{s}_{0}$ $>$ $s_{0}$.

To illustrate the situation, consider the region $R_{+}(\beta _{0})$
surrounded by the $\beta $-axis, the horizontal line $\beta $ $=$ $\beta
_{0},$ and $\Upsilon .$ Observe that $V$ (see (\ref{5.6})) points inward on
the boundary: $\beta $-axis and $\beta $ $=$ $\beta _{0}$ of $R_{+}(\beta
_{0})$ while pointing outward on $\Upsilon $ (see Figure 6.1). The solution $%
p(s)$ moves in $R_{+}(\beta _{0})$ for $\breve{s}_{0}$ $>$ $s$ $>$ $s_{0}.$

\begin{figure}[h]
\includegraphics{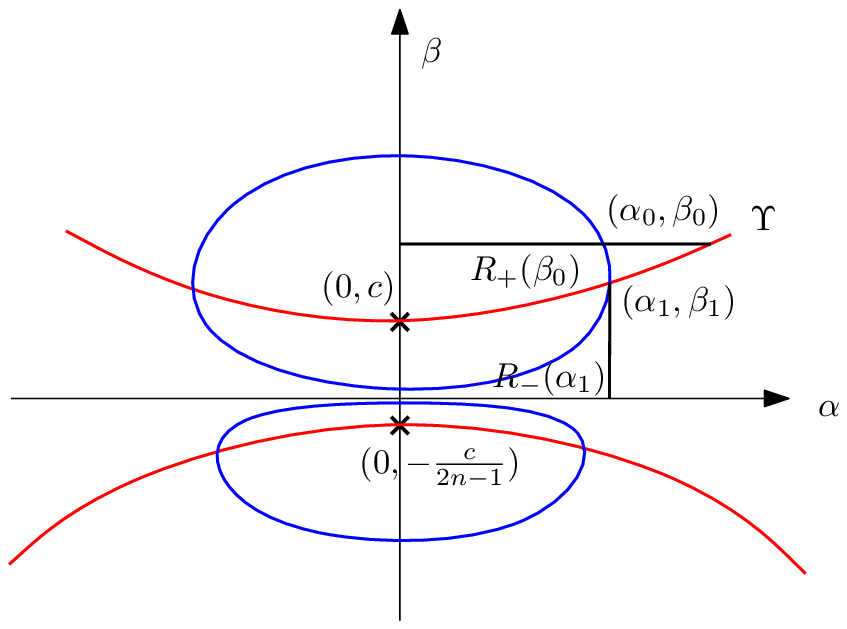}
\caption{Figure 6.1}
\end{figure}

\textbf{Case 3.} $p_{1}$\textit{\ }$=$\textit{\ }$(\alpha _{1},\beta _{1})$ $%
\in $ $\Upsilon \backslash \{(0,c)\}$ $\cap $ (first quadrant). Since $%
\alpha ^{\prime }$ $=$ $0$ and $\beta ^{\prime }$ $<$ $0$ at $p_{1},$ the
solution $p(s)$ to (\ref{5.3}) with $p(s_{1})$ $=$ $p_{1}$ enters $R_{+}$ ($%
R_{-},$ resp.) for a small time interval $s_{1}$ $-$ $\varepsilon <$ $%
s<s_{1} $ ($s_{1}+\varepsilon $ $>$ $s$ $>$ $s_{1},$ resp.$)$.

\textbf{Case 4. }$p_{2}$\textit{\ }$=$\textit{\ }$(\alpha _{2},\beta _{2})$ $%
\in $ $R_{-}.$ Since $\alpha ^{\prime }$ $<$ $0,$ $\beta ^{\prime }$ $<$ $0$
in $R_{-},$ $\alpha (s)$ and $\beta (s)$ are increasing as $s$ changes
towards the negative infinity, where $(\alpha (s),$ $\beta (s))$ $=$ $p(s)$
is the solution to (\ref{5.3}) with $p(s_{2})$ $=$ $p_{2}$. Observe that 
\begin{eqnarray*}
\beta ^{\prime }(s) &=&-2n\beta (s)\alpha (s) \\
&\leq &-2n\beta (s_{2})\alpha (s_{2})<0
\end{eqnarray*}

\noindent for $s$ $\leq $ $s_{2}.$ Suppose $p(s)$ does not hit $\Upsilon $
at any $s$ $<$ $s_{2}.$ Then $\beta (s)$ must go to $+\infty $ as $s$ $%
\rightarrow $ $-\infty .$ So there is $\tilde{s}_{2}$ $<$ $s_{2}$ such that $%
\beta (s)$ $\geq $ $2c$ for $s$ $\leq $ $\tilde{s}_{2}.$ Now from $\alpha
^{\prime }$ $<$ $0$ and (\ref{5.3}), we have%
\begin{eqnarray*}
\alpha ^{2} &\geq &\frac{1}{4n^{2}}(\beta -c)((2n-1)\beta +c) \\
&\geq &\frac{1}{4n^{2}}\frac{\beta }{2}((2n-1)\beta +c) \\
&\geq &\frac{2n-1}{8n^{2}}\beta ^{2}.
\end{eqnarray*}

\noindent It follows that $\alpha $ $\geq $ $c(n)\beta $ where $c(n)$ $=$ $%
\frac{\sqrt{2n-1}}{2\sqrt{2}n}.$ We can then estimate%
\begin{eqnarray*}
-\beta ^{\prime } &=&2n\beta \alpha \\
&\geq &2nc(n)\beta ^{2}=\sqrt{n-\frac{1}{2}}\beta ^{2}
\end{eqnarray*}

\noindent which is reduced to $(\frac{1}{\beta })^{\prime }\geq \sqrt{n-%
\frac{1}{2}}.$ Integrating from $s$ to $\tilde{s}_{2}$ gives%
\begin{equation}
\frac{1}{\beta (\tilde{s}_{2})}-\frac{1}{\beta (s)}\geq \sqrt{n-\frac{1}{2}}(%
\tilde{s}_{2}-s).  \label{5.7}
\end{equation}

\noindent As $s$ $\rightarrow $ $-\infty ,$ the left hand side of (\ref{5.7}%
) is bounded while the right hand side goes to $+\infty .$ The contradiction
shows that $p(s)$ must hit $\Upsilon $ at some finite $\tilde{s}$ $<$ $%
s_{2}. $

On the other hand, consider the region $R_{-}(\alpha _{2})$ surrounded by $%
\Upsilon ,$ $\alpha $ $=$ $\alpha _{2},$ $\beta $ $=$ $0$ ($\alpha $-axis),
and the line segment $\{0\}$ $\times $ $[0,c]$ ($\alpha $ $=$ $0,$ $\ 0$ $%
\leq $ $\beta $ $\leq $ $c)$ (see Figure 6.1)$.$ Observe that the vector
field $V$ (see \ref{5.6}) points inward (towards $R_{-}(\alpha _{2}))$ on $%
\Upsilon ,$ $\alpha $ $=$ $\alpha _{2}$ while pointing outward on $\{0\}$ $%
\times $ $[0,c].$ Note that $V$ does not vanish in $R_{-}(\alpha _{2})$ and $%
\beta $ $=$ $0$ ($\alpha $-axis) is a solution to (\ref{5.3}) with $\alpha
^{\prime }=-\alpha ^{2}+\frac{1}{4n^{2}}(-c^{2})$ $<$ $0.$ Therefore the
solution curve $p(s)$ $:=$ $(\alpha (s),\beta (s))$ to (\ref{5.3}) with $%
p(s_{2})=p_{2}$ must hit either some point in $\{0\}$ $\times $ $(0,c)$ at
finite $\breve{s}_{2}$ $>$ $s_{2}$ or the point $(0,c)$ as $s\rightarrow
+\infty $ by compactness of $\overline{R_{-}(\alpha _{2})}$ and uniqueness
of ($C^{\infty }$ smooth) ODE solutions.

Next suppose $\lim_{s\rightarrow +\infty }p(s)$ $=$ $(0,c)$. We may assume $%
\beta $ $>$ $c$ (otherwise $\beta $ won't tend to $c$ since $\beta $ is
decreasing). From (\ref{5.3}) we compute%
\begin{eqnarray*}
&&\frac{d\alpha ^{\prime }}{ds} \\
&=&-2\alpha \alpha ^{\prime }+\frac{\beta ^{\prime }}{4n^{2}}\{(2n-1)(\beta
-c)+(2n-1)\beta +c\} \\
&\leq &2\alpha ^{3}-\frac{2n\beta \alpha }{4n^{2}}\{(2n-1)\beta +c\}\text{
(by }\beta ^{\prime }<0\text{ and }\beta >c) \\
&\leq &2\alpha ^{3}-c^{2}\alpha .
\end{eqnarray*}

\noindent Since lim$_{s\rightarrow +\infty }\alpha (s)$ $=$ $0,$ we can find
some large number $\breve{s}$ such that $\frac{d\alpha ^{\prime }}{ds}$ $%
\leq $ $0$ for $s$ $\geq $ $\breve{s}$. It follows that $\alpha ^{\prime
}(s) $ $\leq $ $\alpha ^{\prime }(\breve{s})$ $<$ $0.$ But $%
\lim_{s\rightarrow +\infty }\alpha ^{\prime }(s)$ $=$ $0$ ($\alpha $ $%
\rightarrow $ $0,$ $\beta $ $\rightarrow $ $c).$ We have reached a
contradiction. So we conclude that at finite $\breve{s}_{2}$ $>$ $s_{2}$, $p(%
\breve{s}_{2})$ $\in $ $\{0\}$ $\times $ $(0,c).$

\textbf{Case 5}. $p_{+}$ $\in $ $\{0\}$ $\times $ $(0,c).$ Observe that $%
\alpha ^{\prime }$ $<$ $0$ and $\beta ^{\prime }$ $=$ $0$ at $p_{+}.$ The
solution $p(s)$ to (\ref{5.3}) with $p(s_{+})$ $=$ $p_{+}$ will go into the
second quadrant ($R_{-},$ resp.) for a short time after (before, resp.) $%
s_{+}.$

Altogether wherever in the first quadrant we start with, the solution ends
up touching the $\beta $-axis in both finite negative and finite positive
times. Then by the symmetry to the $\beta $-axis we obtain a closed periodic
orbit.

(II) Suppose $p_{0}$\textit{\ }$=$\textit{\ }$(\alpha _{0},\beta _{0})$%
\textit{\ }$\in $\textit{\ }$R^{2,-}\backslash \{(0,-\frac{c}{2n-1})\}.$
Consider the transformation: $\tilde{\alpha}$ $=$ $\alpha ,$ $\tilde{\beta}$ 
$=$ $-\beta .$ Then we have%
\begin{eqnarray}
\tilde{\beta}^{\prime } &=&-2n\tilde{\beta}\tilde{\alpha}  \label{meq} \\
\tilde{\alpha}^{\prime } &=&-\tilde{\alpha}^{2}+\frac{1}{4n^{2}}((2n-1)%
\tilde{\beta}-c)(\tilde{\beta}+c).  \notag
\end{eqnarray}

\noindent Since (\ref{meq}) for $\tilde{\beta}$ $>$ $0$ is similar to (\ref%
{5.3}) for $\beta $ $>$ $0,$ we can analyze (\ref{meq}) similarly to get a
periodic solution $(\tilde{\alpha}(s),$ $\tilde{\beta}(s))$ with $(\alpha
_{0},-\beta _{0})$ as the initial data. Then $(\alpha (s),$ $\beta (s))$ $=$ 
$(\tilde{\alpha}(s),$ $-\tilde{\beta}(s))$ $\in $\textit{\ }$%
R^{2,-}\backslash \{(0,-\frac{c}{2n-1})\}$ is the required periodic
solution. We have completed the proof.

\endproof%

\bigskip

To illustrate the result in Lemma 6.1, please see Figure 6.2 drawn by the
computer.

\begin{figure}[h]
\includegraphics{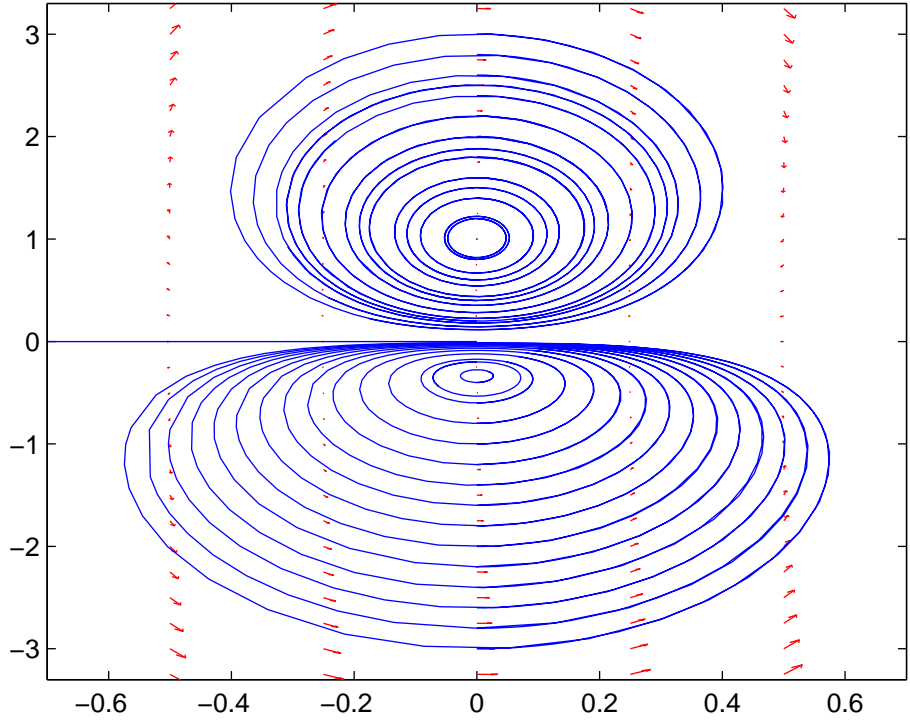}
\caption{Figure 6.2: n=2, c=1}
\end{figure}

\bigskip

\proof
\textbf{(of Lemma B)} Recall that $S_{\Sigma }$ denotes the set of singular
points. For $K$ $\subset $ $\alpha \beta $-plane. we define the subset $%
\Sigma (K)$ $\subset $ $\Sigma \backslash S_{\Sigma }$ by 
\begin{equation*}
\Sigma (K)=\{p\in \Sigma \backslash S_{\Sigma }:(\alpha (p),\beta (p))\in
K\}.
\end{equation*}%
\noindent By Proposition 4.2 we obtain that $k,$ $l$ (and hence $\beta )$,
and $\alpha $ are constant on each leaf of the ($2n-1$)-dimensional
foliation described in Proposition 4.3. On the other hand, $e_{n}$ is
transversal to the leaves by Proposition 4.3, hence $\Sigma (K)$ is open for 
$K$ $=$ $\{(0,c)\}$ or $\{(0,-\frac{c}{2n-1})\}$ or a periodic orbit in the $%
\alpha \beta $-plane, or the $\alpha $-axis by Lemma 6.1. It is also clear
that $\Sigma (K)$ is a closed set for such a $K$. Note that $S_{\Sigma }$
consists of discrete (isolated singular) points by Proposition 4.1. So $%
\Sigma \backslash S_{\Sigma }$ is connected and identified with $\Sigma (K)$
if $\Sigma (K)$ $\neq $ $\emptyset $\ since $\Sigma (K)$ is open and closed.
Observe that $\alpha $ $\rightarrow $ $\pm \infty $ as regular points $p_{j}$
tend to a singular point. For $K$ $=$ $\{(0,c)\}$ or $\{(0,-\frac{c}{2n-1}%
)\} $ or a periodic orbit in the $\alpha \beta $-plane, $\alpha $ is
bounded. Therefore the only choice is $K$ $=$ $\alpha $-axis.(if there
exists a singular point). That is, $0$ $=$ $\beta $ :$=$ $l-2k$ on $\Sigma
\backslash S_{\Sigma }$.

\endproof%

\section{Proof of Theorem A and beyond}

\proof
\textbf{(of Theorem A)} Suppose $\Sigma $ does not contain any singular
point. Then $\Sigma $ is foliated by characteristic curves. Consider the
line field defined by the tangent lines of characteristic curves. Then the
Euler number is the index sum of this line field by Hopf's index theorem (%
\cite{Sp}). Since this line field never vanishes, the Euler number must be
zero. This contradiction to the assumption shows the existence of a singular
point. Next by Lemma B we have $l=2k$ on $\Sigma $. Then by Theorem 1.3, $%
\Sigma $ must be congruent with a Pansu sphere.

\endproof%

\bigskip


Another interesting problem is to relate level sets of a Sobolev extremal to
umbilic hypersurfaces for different Sobolev exponents. Let

\begin{equation*}
p^{\ast }:=\frac{pQ}{Q-p}
\end{equation*}%
\noindent where $Q$ $=$ $2n+2,$ $p$ $\geq $ $1.$ The Sobolev inequality on $%
H_{n}$ reads%
\begin{equation*}
\parallel u\parallel _{L^{p^{\ast }}}\leq C\parallel \nabla _{b}u\parallel
_{L^{p}}
\end{equation*}

\noindent for all functions $u$ such that both sides of the above inequality
are finite. The best constant is obtained by minimizing the Sobolev quotient%
\begin{equation*}
\frac{\parallel \nabla _{b}u\parallel _{L^{p}}}{\parallel u\parallel
_{L^{p^{\ast }}}}.
\end{equation*}

\noindent over all functions $u$ such that both $\parallel \nabla
_{b}u\parallel _{L^{p}}$ and $\parallel u\parallel _{L^{p^{\ast }}}$ are
finite and $\parallel u\parallel _{L^{p^{\ast }}}$ $\neq $ $0.$ The
associated Euler-Lagrange equation reads%
\begin{equation}
\func{div}_{b}(|\nabla _{b}u|^{p-2}\nabla _{b}u)=\sigma u^{p^{\ast }-1}
\label{7.1}
\end{equation}

\noindent where $\sigma $ is a constant (Lagrange multiplier). For other
interesting inequalities on $H_{n},$ the reader is referred to \cite{FL}.

For $p$ $=$ $2,$ equation (\ref{7.1}) is reduced to the $CR$ Yamabe equation%
\begin{equation}
\Delta _{b}u=\sigma u^{1+\frac{2}{n}}.  \label{7.2}
\end{equation}

\noindent Observe that $u(z,t)$ $=$ $(4t^{2}+(|z|^{2}+\lambda )^{2})^{-n/2}$
with constant $\lambda $ $>$ $0$ is a solution to (\ref{7.2}). The level
sets of this solution are "shifted" Heisenberg spheres $\Sigma _{\lambda }$
defined by $4t^{2}$ $+$ $(|z|^{2}+\lambda )^{2}$ $=$ $\rho _{0}^{4}.$
Although these are not Heisenberg spheres (see Example 3.3), they are still
umbilic. Take 
\begin{equation*}
\varphi =\rho _{0}^{4}-[4t^{2}+(|z|^{2}+\lambda )^{2}]
\end{equation*}%
\noindent as a defining function. Let $e_{2n}$ $:=$ $\frac{\nabla
_{b}\varphi }{|\nabla _{b}\varphi |}$ (pointing inwards to the domain $%
\{\varphi $ $>$ $0\}$ at the boundary $\{\varphi $ $=$ $0\})$, $e_{n}$ $:=$ $%
-Je_{2n}$, and $e_{1},$ $..,e_{n-1},$ $e_{n+1},$ $..,$ $e_{2n-1}$ be an
orthonormal frame of $\xi ^{\prime }.$ Then it is not hard to compute $%
h_{jm} $ $=$ $0$ for $1$ $\leq $ $j,m$ $\leq $ $2n-1$ except $j$ $=$ $m$ and 
$|j-m|$ $=$ $n.$ Moreover, we have%
\begin{eqnarray}
l &=&h_{nn}=\frac{2|z|^{2}+(|z|^{2}+\lambda )}{\rho _{0}^{2}|z|},
\label{7.3} \\
k &=&h_{jj}=\frac{|z|^{2}+\lambda }{\rho _{0}^{2}|z|},\text{ }1\leq j\leq
2n-1,j\neq n,  \notag \\
\alpha &=&\frac{2t}{\rho _{0}^{2}|z|}=h_{\beta (n+\beta )}=-h_{(n+\beta
)\beta },\text{ }1\leq \beta \leq n-1.  \notag
\end{eqnarray}

\noindent From (\ref{7.3}) we observe that $l$ $\leq $ $3k$ and 
\begin{equation*}
l=3k\iff \lambda =0\iff \Sigma _{\lambda }\text{ is a Heisenberg sphere.}
\end{equation*}

For $p=1,$ equation (\ref{7.1}) is reduced to the following $p$-mean
curvature equation%
\begin{equation}
H(=\func{div}_{b}\frac{\nabla _{b}u}{|\nabla _{b}u|})=\sigma u^{\frac{1}{2n+1%
}}.  \label{7.4}
\end{equation}

\noindent Observe that $H_{n}\backslash \{0\}$ $=$ $\cup _{0<\lambda <\infty
}$ $S_{\lambda }$ where $S_{\lambda }$ is th Pansu sphere defined in (\ref%
{1.6}). Define a function $u$ on $H_{n}\backslash \{0\}$ by $u$ $=$ $(\frac{%
2n\lambda }{\sigma })^{2n+1}$ on $S_{\lambda }.$ It is not hard to see that $%
u$ $\in $ $C^{2}(H_{n}\backslash \{0\})$ and (\ref{7.4}) holds since, on $%
S_{\lambda },$ $H$ $=$ $2n\lambda $ (see Example 3.2) and $\sigma u^{\frac{1%
}{2n+1}}$ $=$ $2n\lambda $ too. So $u$ is a solution to (\ref{7.4}) with
umbilic level sets $S_{\lambda }.$ In this case, $l$ $=$ $2k.$ We would like
to ask the following question for general $p$ $\geq $ $1:$

\bigskip

\textbf{Question}. \textit{Is each level set of a Sobolev extremal, solution
to (\ref{7.1}), umbilic?}


\bigskip

\begin{thebibliography}{99}
\bibitem{Al} Alexandrov, A. D., \textit{Uniqueness theorems for surfaces in
the large I}, Vestnik Leningrad Univ., 11 (1956) 5-17.

\bibitem{CH04} Cheng, J.-H. and Hwang, J.-F., \textit{Properly embedded and
immersed minimal surfaces in the Heisenberg group}, Bull. Aus. Math. Soc.,
70 (2004) 507-520.

\bibitem{CH10} Cheng, J.-H. and Hwang, J.-F., \textit{Variations of
generalized area functionals and p-area minimizers of bounded variation in
the Heisenberg group}, Bulletin of the Institute of Mathematics, Academia
Sinica, New Series, 5 (2010) 369-412.

\bibitem{CH14} Cheng, J.-H. and Hwang, J.-F., \textit{Uniqueness of
generalized }$p$\textit{-area minimizers and integrability of a horizontal
normal in the Heisenberg group}, Calc. Var. and PDE;
http://arxiv.org/abs/1211.1474 (published online 2013).

\bibitem{CHMY04} Cheng,~J.-H.,~Hwang,~J.-F.,~Malchiodi,~A.,~and~Yang,~P., 
\textit{Minimal surfaces in pseudohermitian geometry}, Annali della Scuola
Normale Superiore di Pisa, Classe di Scienze~(5), 4 (2005) 129-177.

\bibitem{CHMY12} Cheng,~J.-H.,~Hwang,~J.-F.,~Malchiodi,~A.,~and~Yang,~P., 
\textit{A Codazzi-like equation and the singular set for }$C^{1}$\textit{\
smooth surfaces in the Heisenberg group}, Journal fur die reine und
angewandte Mathematik, 671 (2012) 131-198.

\bibitem{CHY1} Cheng, J.-H., Hwang, J.-F., and Yang, P., \textit{Existence
and uniqueness for }$p$\textit{-area minimizers in the Heisenberg group},
Math. Annalen, 337 (2007) 253-293.

\bibitem{CHY2} Cheng, J.-H., Hwang, J.-F., and Yang, P., \textit{Regularity
of }$C^{1}$\textit{\ smooth surfaces with prescribed }$p$\textit{-mean
curvature in the Heisenberg group}, Math. Annalen, 344 (2009) 1-35.

\bibitem{CL} Chiu,~H.-L.~and~Lai,~S.-H., \textit{The fundamental theorem for
hypersurfaces in Heisenberg groups}, Cal. Var. and P.D.E., DOI
10.1007/s00526-015-0818-1, 2015.

\bibitem{FL} Frank, R. L. and Lieb, E. H., \textit{Sharp constants in
several inequalities on the Heisenberg group}, Ann. of Math., 176 (2012)
349-381.

\bibitem{Lee} Lee,~J.~M., \textit{The Fefferman metric and pseudohermitian
invariants}, Trans. Amer. Math. Soc., 296 (1986) 411-429.

\bibitem{LM} Lin, Y., and Ma, H., \textit{A characterization of spheres in
the Heisenberg group }$H^{n},$ preprint.

\bibitem{MR} Montiel, S. and Ros, A., \textit{Compact hypersurfaces: the
Alexandrov theorem for higher order mean curvatures}, Diff. Geom., A
symposium in honour of Manfredo do Carmo, Pitman Mono. and Surv. in Pure and
Appl. Math. 52, ed. B. Lawson and K.Tenenblat, pp 279-296.

\bibitem{Pau04} Pauls, S. D., \textit{Minimal surfaces in the Heisenberg
group}, Geometric Dedicata, 104 (2004) 201-231.

\bibitem{Pauls} Pauls, S. D., \textit{H-minimal graphs of low regularity in }%
$H^{1},$ Comment. Math. Helv. 81 (2006) 337-381; arXiv: math.DG/0505287 v3,
Nov. 1, 2006 (to which the reader is referred).

\bibitem{Ri} Ritor\'{e}, M., \textit{Examples of area-minimizing surfaces in
the subriemannian Heisenberg group }$H^{1}$\textit{\ with low regularity},
Calc. Var. and PDE (2008), doi:10.1007/s00526-008-0181-6

\bibitem{Rit} Ritor\'{e},~M., \textit{A proof by calibration of an
isoperimetric inequality in the Heisenberg group }$\mathit{H}^{n}$, Calc.
Var. and PDE, 44 (2012) 47-60.

\bibitem{RR} Ritor\'{e}, M. and Rosales, C., \textit{Area-stationary
surfaces in the Heisenberg group }$H^{1},$ Advances in Math., 219 (2008)
633-671.

\bibitem{Sp} Spivak,~M.,\textit{\ A comprehensive introduction to
differential~geometry}, Vol. 3, Publish or Perish, Inc., Boston,~1975.

\bibitem{Web} Webster, S. M., \textit{Pseudohermitian structures on a real
hypersurface}, J. Diff. Geom. 13 (1978) 25-41.
\end{thebibliography}
\end{document}